\documentclass[10pt]{article} 

\usepackage{geometry}
\usepackage{booktabs} 
\usepackage{array}
\usepackage{subcaption}
\usepackage{amsmath,amssymb,amsthm,bm,bbm}
\usepackage[bb=esstix]{mathalpha}
\usepackage{mathtools}
\usepackage{natbib}
\usepackage{epsfig,graphicx}
\usepackage{comment}
\usepackage{color}
\definecolor{blue(pigment)}{rgb}{0.2, 0.2, 0.6}
\definecolor{ultramarine}{rgb}{0.07, 0.04, 0.56}
\definecolor{darkspringgreen}{rgb}{0.09, 0.45, 0.27}
\definecolor{hookersgreen}{rgb}{0.0, 0.44, 0.0}
\definecolor{plum(traditional)}{rgb}{0.56, 0.27, 0.52}
\definecolor{purple(html/css)}{rgb}{0.5, 0.0, 0.5}
\definecolor{magenta(dye)}{rgb}{0.79, 0.08, 0.48}
\usepackage{srcltx}
\usepackage[mathscr]{eucal}
\usepackage[math]{easyeqn}
\usepackage{etoolbox}
\usepackage{hyperref}
\hypersetup{
            colorlinks,
            linkcolor=hookersgreen,
            linktoc=blue,
            citecolor=ultramarine,
            urlcolor=black,
            filecolor=black
            }
\usepackage{nameref}
\usepackage{multirow}
\usepackage{accents}
\usepackage{mathrsfs}
\usepackage{algorithm}
\usepackage{algpseudocode}
\DeclareMathAlphabet{\mathbbmsl}{U}{bbm}{bx}{sl}
\newcommand{\bb}[1]{\boldsymbol{#1}}
\def\argmin{\operatornamewithlimits{argmin}}
\def\eqdef{\stackrel{\operatorname{def}}{=}}
\renewcommand{\hat}[1]{\widehat{#1}}
\def\uv{\bb{u}}
\def\vv{\bb{v}}
\def\wv{\bb{w}}
\def\xv{\bb{x}}
\def\yv{\bb{y}}
\def\ev{\bb{e}}
\def\l{\left}
\def\r{\right}

\def\lgd{f}

\def\PfL{\P_{\lgd}}

\def\vH{\vA}

\def\dltw{\delta}

\def\dltwb{\omega}

\def\dltwu{\tau}

\def\vv{\bb{v}}

\def\DFL{\mathbb{D}}
\def\DFLG{\DFL_{\GP}}

\def\R{\mathbbmsl{R}}
\def\E{\mathbbmsl{E}}

\def\P{\mathbbmsl{P}}

\def\kappa{\varkappa}

\def\diag{\operatorname{diag}}

\def\ND{\mathcal{N}}

\def\T{\top}
\def\tr{\operatorname{tr}}
\def\Tr{\operatorname{Tr}}
\def\TV{\operatorname{TV}}

\def\ex{\mathrm{e}}

\def\Id{\mathbbmsl{I}}
\def\Ind{\operatorname{1}\hspace{-4.3pt}\operatorname{I}}

\def\BBH{B}

\def\DFc{\DF_{0}}

\def\dimA{\mathtt{p}}

\def\dimL{\dimA}

\def\dimH{\dimA}

\def\dimp{p}

\def\gaussv{\bb{\gauss}}

\def\GP{G}

\def\normG{\alpha}
\def\rr{\mathtt{r}}

\def\rrL{\rr}

\def\rrdL{\rrL}

\def\Tau{T}

\def\UV{\mathcal{U}}

\def\vA{\mathtt{v}}

\def\wv{\bb{w}}

\def\xx{\mathtt{x}}

\def\zq{z}

\newcommand{\beq}{\begin{equation}}
\newcommand{\eeq}{\end{equation}}
\def\beqs#1\eeqs{
    \begin{equation}\begin{split}
    #1
    \end{split}\end{equation}
}
\def\beqsj #1\eeqsj
{
    \begin{equation}
    \setlength{\jot}{10pt}
    \begin{split}
    #1
    \end{split}\end{equation}
}

\def\beqsn#1\eeqsn{
    \begin{equation*}\begin{split}
    #1
    \end{split}\end{equation*}
}

\newcommand{\e}{\exp}

\newcommand{\les}{\lesssim}

\renewcommand{\l}{\left}
\renewcommand{\r}{\right}

\renewcommand{\gaussv}{Z}

\renewcommand{\DFc}{\DFL}
\newcommand{\rrdc}{\rrdL}

\renewcommand{\Tau}{D}

\renewcommand{\gaussv}{Z}
\renewcommand{\Tau}{D}

\renewcommand{\gaussv}{Z}

\renewcommand{\PfL}{\pi_{\lgd}}
\newcommand{\lap}{\gamma_{\lgd}}

\renewcommand{\lgd}{f}
\renewcommand{\PfL}{\pi_\lgd}
\newcommand{\map}{\hat{\thv}}
\newcommand{\thv}{\bm\theta}
\newcommand{\calR}{\mathcal R}
\newcommand{\calS}{\mathcal S}

\numberwithin{equation}{section}
\numberwithin{figure}{section}
\newcounter{example}[section]
\numberwithin{example}{section}
\newcounter{remark}[section]
\numberwithin{remark}{section}
\newtheorem{theorem}{Theorem}[section]
\newtheorem{proposition}[theorem]{Proposition}
\newtheorem{lemma}[theorem]{Lemma}
\newtheorem{corollary}[theorem]{Corollary}
\newtheorem{definition}[theorem]{Definition}
\newtheorem{exmp}[example]{Example}
\newtheorem{rmrk}[remark]{Remark}
\newenvironment{example}{\begin{exmp}\rm}{\end{exmp}}
\newenvironment{remark}{\begin{rmrk}\rm}{\end{rmrk}}
\newtheorem{assumption}[theorem]{Assumption}

\title{A unified theory of the high-dimensional Laplace approximation with application to Bayesian inverse problems}

\author{%
  \begin{minipage}{0.45\textwidth}
    \centering
    Anya Katsevich\thanks{This work was partially supported by NSF grant DMS-2202963.} \\
    \textit{Department of Statistical Science} \\
    \textit{Duke University} \\
    \texttt{anya.katsevich@duke.edu}
  \end{minipage}%
  \hfill%
  \begin{minipage}{0.45\textwidth}
    \centering
    Vladimir Spokoiny\thanks{Financial support by the German Research Foundation (DFG) through the Collaborative Research Center 1294 ``Data assimilation'' is gratefully acknowledged.} \\
    \textit{Weierstrass Institute and HU Berlin} \\
    \texttt{spokoiny@wias-berlin.de}
  \end{minipage}
}

\date{\today}

\begin{document}

\maketitle

\begin{abstract}
The Laplace approximation (LA) to posteriors is a ubiquitous tool to simplify Bayesian computation, particularly in the high-dimensional settings arising in Bayesian inverse problems. Precisely quantifying the LA accuracy is a challenging problem in the high-dimensional regime. We develop a theory of the LA accuracy to high-dimensional posteriors which both subsumes and unifies a number of results in the literature. The primary advantage of our theory is that we introduce a new degree of flexibility, which can be used to obtain problem-specific upper bounds which are much tighter than previous ``rigid" bounds. We demonstrate the theory in a prototypical example of a Bayesian inverse problem, in which this flexibility enables us to improve on prior bounds by an order of magnitude. Our optimized bounds in this setting are dimension-free, and therefore valid in arbitrarily high dimensions.\end{abstract}

\noindent
\textbf{AMS 2020 Subject Classification:} 62F15, 62E17

\vspace{0.5em}

\noindent
\textbf{Keywords:} effective dimension, Laplace approximation, Bayesian inverse problems

\section{Introduction}

The Bayesian approach to inverse problems offers a natural and convenient framework to do uncertainty quantification. However, due to the high-dimensional nature of most inverse problems, the Bayesian approach can be very computationally intensive~\cite{kaipio2006statistical, stuart2010inverse, dashti2017bayesian}. A standard way to simplify the computations involved is to use the Laplace approximation (LA), a Gaussian approximation to the posterior first introduced in the Bayesian context by~\cite{tierney1986accurate}. The LA has been used to great effect in Bayesian inverse problems even at extreme scale~\cite{stadler2012extreme,ghattas2021learning}. But to rely on the LA for uncertainty quantification, its approximation accuracy must be ascertained. Determining the LA accuracy to the posterior in the most relevant, high-dimensional regime, is a challenging theoretical endeavor and currently an active area of research.

There have been a number of works addressing this question, including~\cite{huggins2018practical,dehaene2019deterministic,schillings2020convergence,helin2022non,fischer2022normal}. The strongest guarantees have been obtained by~\cite{bp},~\cite{katsBVM}, and~\cite{spok23}. The work~\cite{katskew} obtains similar guarantees but under stronger assumptions. We postpone discussion of this latter work to Section~\ref{related}.

To describe the prior work and our own contributions in more detail, we specify the problem mathematically. Consider a posterior probability density of the form \(\PfL(\thv)\propto e^{-\lgd(\thv)}\) on \(\R^{\dimp}\). Suppose \(\lgd\) has a unique global minimizer \(\map\). Then the LA \(\lap\) to \(\PfL\) is defined as the Gaussian distribution
\begin{EQA}
&&\lap=\ND(\map,\,\nabla^2\lgd(\map)^{-1}).
\end{EQA} 
The function \(\lgd\) depends on a large parameter \(n\), which denotes either sample size, or inverse noise level, or number of measurements.  Roughly speaking, the Hessian \(\nabla^2\lgd\) grows linearly with \(n\). 

We now give a unified presentation of the bounds of the three works~\cite{bp},~\cite{katsBVM},~\cite{spok23}, on the total variation (TV) distance between \(\PfL\) and \(\lap\). The form in which we present the first two bounds is different from how they are presented in those works, for reasons which will soon become clear. In fact, this new presentation of the bounds is one of the insights of our work.


~\cite{katsBVM} shows that \(\TV(\PfL,\lap)\lesssim \dltwu_3\dimp\). Here, recall that \(\dimp\) is the dimension of the unknown parameter \(\theta\). The quantity \(\dltwu_3\) depends on the operator norm of \(\nabla^3\lgd\). The work~\cite{bp} proves \(\TV(\PfL,\lap)\lesssim \dltwu_3'\dimp'\). The definition of \(\dltwu_3'\) is similar but not identical to \(\dltwu_3\), and \(\dimp'\) can be thought of as an \emph{effective} dimension, though this terminology is not used in~\cite{bp}. Finally,~\cite{spok23} shows \(\TV(\PfL,\lap)\lesssim \dltwu_3''{\dimp''}^{3/2}\), where \(\dltwu_3''\) is analogous but not identical to \(\dltwu_3,\dltwu_3'\). The term \(\dimp''\) is called the ``Laplace effective dimension".  Both \(\dimp'\) and \(\dimp''\) are always less than or equal to the full dimension \(\dimp\), and~\cite{spok23} notes that \(\dimp''\) can in some cases be significantly smaller than \(\dimp\). The terms \(\dltwu_3,\dltwu_3',\dltwu_3''\) all scale with \(n\) as \(n^{-1/2}\), though they are also sensitive to the statistical model under consideration, and may grow with dimension or effective dimension. 

These three bounds are the strongest in the literature. Indeed, prior works required \(n\gg\dimp^3\) for accuracy of the LA (e.g.~\cite{helin2022non}),  whereas~\cite{katsBVM,bp} only require \(n\gg\dimp^2\) or \(n\gg{\dimp'}^2\), and~\cite{spok23} only requires \(n\gg {\dimp''}^3\) (disregarding model-dependent terms). However, there has been no clear understanding of how the bounds relate to one another. For example, although \(\dimp'\leq\dimp\), the bound of~\cite{bp} is not strictly tighter than that of~\cite{katsBVM}, because it is not always true that \(\dltwu_3'\leq \dltwu_3\). And if \(\dimp''\) is potentially much smaller than \(\dimp\), then the bound of~\cite{spok23} may be tighter than that of~\cite{katsBVM} despite the higher power of \(\dimp''\). 

Yet the fact that we are able to write all three bounds in a very similar form strongly suggests that there is a deep connection between them. In this paper, we show that this is indeed the case. We develop a unified theory from which the three bounds follow, which explains the relationship among them, and which allows us to obtain more powerful estimates than each of these.  
We prove that, for all \(\DFL\succ0\) in a certain broad class of matrices satisfying \(\|\nabla^2\lgd(\map)^{-1/2}\DFL\|=1\), it holds 
\begingroup
 \setlength{\jot}{10pt}
 \begin{equation}\label{TV1}\begin{gathered}
 \TV(\PfL,\lap)\lesssim \dltwu_3(\DFL)\dimL(\DFL),\\
\dltwu_3(\DFL)\eqdef\sup_{\thv\in\UV}\|\nabla^3\lgd(\thv)\|_{\DFL},\qquad \dimL(\DFL)\eqdef\Tr(\nabla^2\lgd(\map)^{-1}\DFL^2).
 \end{gathered}\end{equation} \endgroup The set \(\UV\) is a neighborhood of the minimizer \(\map\) on which \(\PfL\) and \(\lap\) both concentrate, and \(\|\nabla^3\lgd(\thv)\|_{\DFL}\) is the operator norm with respect to the Mahalanobis norm \(\|\DFL\cdot\|\) on \(\R^{\dimp}\). The term \(\dltwu_3(\DFL)\) expresses how far the log posterior is from quadratic, and scales with the large parameter \(n\) as \(n^{-1/2}\). The term \(\dimL(\DFL)\) is a notion of effective dimension of the posterior (note that \(1\leq\dimL(\DFL)\leq\dimp\)). In inverse problem settings and for prudent choices of \(\DFL\), it captures the number of coordinates for which the information in the likelihood outweighs the prior regularization strength. This number is typically significantly smaller than the full parameter dimension \(\dimp\). Up to model-dependent factors affecting the magnitude of \(\dltwu_3(\DFL)\), the bound~\eqref{TV1} is small provided \(n\gg\dimL(\DFL)^2\). 
  
Furthermore, we show that \emph{the bounds of both~\cite{katsBVM} and~\cite{bp} can be brought into the exact form}~\eqref{TV1} but with two particular matrices \(\DFL\). The bound of~\cite{spok23} is of this same form but with \(\dimL(\DFL)^{3/2}\) instead of \(\dimL(\DFL)\), and yet a third matrix \(\DFL\). Thus, the quantities \(\dltwu_3, \dltwu_3',\dltwu_3''\) described above correspond to \(\dltwu_3(\DFL)\) for three fixed choices of \(\DFL\), and similarly for \(\dimp,\dimp',\dimp''\), with the same respective \(\DFL\). In~\cite{spok23}, the ``default" choice is to take \(\DFL^2\) to be the Hessian of the negative log likelihood at the point \(\map\), but the formulation of the results does technically allow for the possibility of choosing \(\DFL\) differently.

The main advantage of our theory is the introduction of the free parameter \(\DFL\) into the bound, which can be optimized. There is a trade-off between the two factors: increasing \(\DFL\) in the positive definite ordering decreases \(\dltwu_3(\DFL)\) and increases \(\dimL(\DFL)\), so an optimal middle ground can be found. Our work reveals that the bounds of~\cite{katsBVM} and~\cite{bp} are simply two choices within this spectrum. They are not directly comparable to one another, but are instead both subsumed by our result~\eqref{TV1}. 

We demonstrate the power of our new, more flexible approach in the setting of an inverse problem. We choose an optimal \(\DFL\) which leads to a TV bound that is tighter by an order of magnitude than those of~\cite{katsBVM} and~\cite{bp}. We don't compare our results to~\cite{spok23} due to the inherently worse dimension dependence.


 In summary, our main contributions are the following.
 \begin{enumerate}
 \item We develop a theory of the LA accuracy which subsumes and unifies the best estimates to date.
 \item The power of the theory is that it can be tailored to the specific problem at hand via the choice of \(\DFL\).
 \item We demonstrate the theory in a prototypical example of an inverse problem, in which this flexibility enables us to improve on prior bounds by an order of magnitude. 
 \end{enumerate}

\subsection{Application to an inverse problem}\label{intro:inv}
Informally, the problem we study is to recover the first \(\dimp\) coefficients of an orthonormal basis representation of a function \(q^*\) given noisy perturbations of \((\calR q^*)(j/n)\), \(j=1,\dots,n\). The operator \(\calR\) is a linear smoothing operator, so that recovering \(q^*\) from \(\calR q^*\) is an inverse problem. More precisely, we consider the following generalized linear inverse problem: we are given independent data 
 \begin{EQA}\label{Yjintro}
&& Y_j \sim \rho(\cdot \mid (\calR q^*)(j/n)),\qquad j=1,\dots,n,
 \end{EQA} where \(\rho(\cdot\mid\cdot)\) is an exponential family and \(\calR\) is a linear operator from \(L^2[0,1]\) to itself, such that the eigenvalues of \(\calR^\T\calR\) decay as \(k^{-2\beta}\), \(k=1,2,3,\dots\). We wish to recover the function \(q^*\), which we approximate in the basis of eigenfunctions of \(\calR^\T\calR\), using \(\dimp\) coefficients. Finally, we put a diagonal Gaussian prior  \(\ND(0,G^{-2})\) over these coefficients, with \(G^2=\diag(k^{2\gamma})_{k=1}^{\dimp}\). Unless \(\rho\) is Gaussian, we obtain a non-Gaussian posterior over the \(\dimp\) coefficients. Thus the LA, a Gaussian approximation, can be used to simplify computation with this posterior. Our goal is to determine the LA accuracy in this setting.

 Our model is inspired by imaging modalities such as positron emission tomography (PET), which are based on X-rays (high-energy photons). Our \(\calR q^*\) is analogous to line integrals of the intensity distribution \(q^*\) of radiation sources, and the arguments \(j/n\) correspond to different lines. Typically, the measurement error of these line integrals  (our \(\rho\)) obeys a Poisson distribution~\cite{hohage2016inverse}. See this work for an overview of inverse problems with Poisson noise. The work~\cite{melidonis2023efficient} notes that in addition to Poisson noise, geometric and binomial noise are also appropriate modeling assumptions for certain low-photon imaging problems. These also give rise to valid choices of \(\rho\) in our setting.
 



We compute the LA TV error bound as a function of the problem parameters, both for our optimized choice of \(\DFL\) and for the choices of \(\DFL\) from prior work. Our optimized bound improves on the bounds from prior works by an order of magnitude. Furthermore, under mild assumptions on the smoothing and regularization parameters \(\beta,\gamma\), the bound can be expressed purely as some negative power of \(n\), so in particular, it is dimension free. Thus the bound is valid for arbitrarily large \(\dimp\).

As has been discussed in~\cite{katskew} and~\cite{bp}, LA error bounds are typically challenging to compute due to the presence of tensor operator norms. Our application is no exception. In fact, even for the standard choices of \(\DFL\) used in prior works, the term \(\dltwu_3(\DFL)\) is in general intractable in our setting.  The quantity \(\dimL(\DFL)\) is also challenging to compute, particularly when \(\nabla^2\lgd(\map)\) and \(\DFL^2\) are not simultaneously diagonalizable. To handle these difficulties, we impose an additional assumption on \(\calR\) which guarantees that \(\nabla^2\lgd(\map)\) and \(\DFL^2\) are both close to diagonal matrices (for all considered choices of \(\DFL\)). This significantly simplifies the computation of \(\dltwu_3(\DFL)\) and \(\dimL(\DFL)\). One of our main technical achievements is to show this assumption on \(\calR\) is verified for a large class of integral operators. This requires studying the asymptotics of eigenfunctions of \(\calR\calR^\T\), and we use Sturm-Liouville theory to do so. 



To conclude, we highlight that our LA error analysis in this inverse problem setting is important in its own right, and not just as a demonstration of the general theory. We have elucidated how the LA accuracy depends on the key problem parameters: \(n\) (number of measurements of \(\calR q^*\)), \(\dimp\) (number of basis functions used to represent \(q^*\)), \(\beta\) (smoothness of \(\calR\)), and \(\gamma\) (regularization strength) for an important class of Bayesian inverse problems. Although we have made a number of simplifying assumptions, including that the regularization is in the basis of eigenfunctions of \(\calR^\T\calR\), our analysis is a crucial first step toward understanding more realistic settings.

\subsection{Related work}\label{related}
The work~\cite{katskew} derives a leading order decomposition of the LA TV error, of the form \(\mathrm{TV}(\PfL,\lap)=L_{\mathrm{TV}} + R_{\mathrm{TV}}\). The term \(L_{\mathrm{TV}}\) only depends on \(\nabla^2\lgd(\map)\) and \(\nabla^3\lgd(\map)\), and the term \(R_{\mathrm{TV}}\) is bounded in terms of the second, third, and fourth derivatives.  The decomposition always yields an upper bound on \(\mathrm{TV}(\PfL,\lap)\), and in some cases also provides lower bounds, when one can show that \(R_{\mathrm{TV}}\ll L_{\mathrm{TV}}\). The upper bound belongs to the group of recent, strongest results, in the sense that it generically requires \(n\gg\dimp^2\), rather than \(n\gg\dimp^3\). However, since it derives from an expansion, the structure of the bound is different from the structure~\eqref{TV1} common to the present work and that of~\cite{katsBVM,bp}. In particular, using the bound requires computing the fourth derivative. Furthermore, in the inverse problem setting described above, the bound of~\cite{katskew} is not dimension-free, unlike our bound with an optimized choice of \(\DFL\); see Remark~\ref{rk:katskew}. Because of these differences, we do not provide a detailed analysis of the bound of~\cite{katskew} in the present work.

 The LA is not the only way to approximate a posterior by a simpler distribution. For example, the works~\cite{nickl2022polynomial,bohr2024concave} prove bounds on the accuracy of log-concave approximations to posteriors for non-linear Bayesian inverse problems.

Frequentist properties of posteriors in Bayesian inverse problems have been established in numerous works. For Bayesian linear inverse problems,~\cite{knapik2011bayesian,knapik2018general} study posterior contraction around the true parameter, as well as frequentist coverage of credible sets, in infinite dimensions.~\cite{bochkina2013consistency} studies posterior contraction around the true parameter for the same model studied here, that of a generalized linear inverse problem.

Another type of result related to LA error bounds is the Bernstein von-Mises (BvM) theorem, which also shows the posterior is asymptotically Gaussian, but from a frequentist perspective.~\cite{monard2019efficient,nickl2020bernstein} prove the BvM for linear and nonlinear inverse problems in infinite dimensions, respectively.~\cite{lu2017bernstein} proves the BvM for nonlinear inverse problems in high but finite dimensions. 
See~\cite{nickl2023bayesian} for an in-depth treatment of the statistical theory and computational aspects of the Bayesian approach to nonlinear inverse problems.

We note that in contrast, we do not take the frequentist perspective in this work. The object of interest is the posterior itself and the accuracy of the LA to the posterior, rather than an underlying ground truth parameter. 
\subsection*{Organization}
In Section~\ref{sec:results}, we present our main results in the general setting, and compare them to prior work. In addition to the aforementioned TV bound, another one of our results is on posterior concentration. In Section~\ref{inv} we present our results on the LA accuracy for the generalized linear inverse model, and compare our bound to those from prior works, when applied to this model. 
Omitted proofs can be found in the Appendix.

\subsection*{Notation}Let \(T\) be a \(\dimp\times\dimp\times\dimp\) tensor and \(\DFL\) a \(\dimp\times\dimp\) positive definite matrix. We let
\begin{EQA}\label{TDdef}
\|T\|_{\DFL} &=& \sup_{\|\DFL\uv\|=\|\DFL\vv\|=\|\DFL\wv\|=1}\langle T, \uv\otimes\vv \otimes\wv\rangle.
\end{EQA} 

\section{Main results in the general setting}\label{sec:results}

In Section~\ref{SsetupLapl}, we describe the general problem setting. We give our result on posterior concentration in Section~\ref{Sconc}, and on the LA TV error in Section~\ref{STV}. In Section~\ref{sec:comparegen} we compare our results to those of~\cite{spok23,katsBVM,bp}. Finally, we present the proof of our TV bound in Section~\ref{sec:TVproof}.

\subsection{Set-up}\label{SsetupLapl}


We consider a probability density 
\begin{EQA}
&&\PfL(\thv)\propto \ex^{-\lgd(\thv)},\qquad\thv\in\R^{\dimp},
\end{EQA} where \(\dimp\) may be very large but is not infinite. We assume \(\PfL\) has a unique MAP (maximum aposteriori) \(\map\). Equivalently, \(\map\) is the unique global minimizer of  \(\lgd\):
\begin{EQA}
	\map
	&=&
	\argmin_{\thv \in \R^{\dimp}} \lgd(\thv) .
\label{scdygw7ytd7wqqsquuqydtdtd}
\end{EQA}
We also assume \( \lgd \) is at least three times differentiable. Let 
\begin{EQA}
&&\DFLG^2=\nabla^2\lgd(\map),
\end{EQA} the Hessian of \(\lgd\) at its minimum \(\map\), where \(\DFLG\succ0\) is symmetric. The \emph{Laplace approximation} (LA) to \(\PfL\) is the Gaussian distribution \(\ND(\map, \DFLG^{-2})\), with probability density
\begin{EQA}\label{df:lap}
\lap(\thv) &\propto&  \exp(-\|\DFLG(\thv-\map)\|^2/2).
\end{EQA}
Comparing various matrices \(\DFL^2\succ0\) to the Hessian \(\DFLG^2\) plays a central role in our analysis. We introduce two very important quantities to this effect:
\begin{definition}[Norm and ``effective dimension" of \(\DFL\)]\label{def:effdim} Let \(\DFL\succ0\) be a symmetric matrix. We define
\begin{EQA}
&&\normG(\DFL)=\|\DFLG^{-1}\DFL\| = \|\DFL\DFLG^{-1}\|=\sqrt{\|\DFLG^{-1}\DFL^2\DFLG^{-1}\|},
\end{EQA} and 
\begin{EQA}
&&\dimL(\DFL) = \frac{\Tr(\DFLG^{-1}\DFL^2\DFLG^{-1})}{\|\DFLG^{-1}\DFL^2\DFLG^{-1}\|} =  \frac{\Tr(\DFLG^{-2}\DFL^2)}{\|\DFLG^{-1}\DFL\|^2} =\Tr(\DFLG^{-2}(\DFL/\normG(\DFL))^2).\label{eq:effdim}
\end{EQA}
\end{definition}
\noindent Note that \(1\leq \dimL(\DFL)\leq\dimp\), and that \(\dimL(\DFL)=\dimL(c\DFL)\) for any \(c>0\). Also, we have \(\dimL(\DFLG)=\dimp\).  A key component of our approach will be to choose \(\DFL\neq \DFLG\) to ensure that \(\dimL(\DFL)\ll \dimp\). We give two useful characterizations of \(\normG(\DFL)\):
\begin{EQA}\label{normGequiv}
&&\normG(\DFL)^2 = \inf\{\lambda>0\, :  \,\DFL^2\preceq \lambda\DFLG^2\},\\
&&\normG(\DFL)^2=\lambda\quad\Longleftrightarrow \quad \DFL^2\preceq \lambda\DFLG^2,\;\exists\uv\neq0, \|\DFL\uv\|^2=\lambda\|\DFLG\uv\|^2.
\end{EQA}
\begin{remark}Our definition of effective dimension differs from that of~\cite{spok23} in that here, we normalize by \(\|\DFLG^{-1}\DFL\|^2\). A conceptual difference between the two works is that here, we think of effective dimension not as an ``absolute" quantity, but one that depends on \(\DFL\). Nevertheless, as noted in Remark~\ref{rk:effd}, we observe that \(\dimL(\DFL)\) is very stable to changes in \(\DFL\) in the inverse problem setting. 
\end{remark}

\subsubsection{Structure of \(\lgd\)}\label{subsub:struct}
To conclude this section, we describe the structure of \(\lgd\) satisfied in our main application of a generalized linear inverse model. First, in the general context of Bayesian inference, the posterior takes the form \(\PfL(\thv)\propto P(Y^n\mid\thv)\pi_0(\thv)\), where \(P(Y^n\mid\thv)\) is the likelihood for observed data \(Y^n\), and \(\pi_0\) is the prior probability density. We let \(L\) denote the negative log likelihood: \(L(\thv)=-\log P(Y^n\mid\thv)\). Furthermore, we primarily consider a centered Gaussian prior \(\mathcal N(0, G^{-2})\) for some symmetric \(G\succ0\), in which case \(\pi_0(\thv)\propto\exp(-\|G\thv\|^2/2)\). This leads to the following decomposition of the function \(\lgd\):
\begin{EQA}\label{fLG}
&&\lgd(\thv)=L(\thv)+\frac12\|G\thv\|^2=-\log\PfL(\thv)+\mathrm{const.}
\end{EQA} In the setting arising in this work, the function \(L\) is convex but not strongly convex. Thus \(\nabla^2L\succeq0\) on \(\R^{\dimp}\), but there is no uniform positive definite lower bound on the Hessian of \(L\). We note that the Hessian \(\DFLG^2\) of \(\lgd\) at the minimum is then given by
\begin{EQA}\label{DFLG-L-G}
&&\DFLG^2=\nabla^2L(\map)+G^2.
\end{EQA}
Finally, informally, the negative log likelihood \(L\) scales linearly with sample size \(n\). However, this statement should be interpreted with care, because in the inverse problem context, the dependence of \(L\) on \(n\) is anisotropic. For example, in the generalized linear inverse model considered in Section~\ref{inv}, we show that 
\begin{EQA}\label{hess-L-n}
&&\nabla^2L(\map) \asymp \diag\left(\left\{\frac{n}{k^{2\beta}}\right\}_{k=1}^{\dimp}\right).
\end{EQA} If \(\dimp\) grows large with \(n\), then for large values of the index \(k\), it is not true that \(n/k^{2\beta}\geq Cn\) as \(n\to\infty\). This is in contrast to other more benign statistical settings outside of the inverse problem context. 

Note that~\eqref{fLG} and convexity of \(L\) imply \(\lgd\) is strongly convex. Although we use strong convexity in our concentration result in Section~\ref{Sconc}, strong convexity is \emph{not} necessary for our TV bound in Section~\ref{STV}.

\subsection{Concentration}\label{Sconc}
In this section, we are interested in determining the radius $\rrdc$ such that $\PfL$ concentrates in the set
\begin{EQA}\label{def:UV-DFL}
&&\UV(\DFc,\rrdc) = \{\thv\in\R^{\dimp}\; :\|\DFc(\thv-\map)\|\leq\rrdc\}.
\end{EQA} Here, $\DFc\succ0$ is some user-specified symmetric matrix, and the radius at which $\UV(\DFc,\rrdc)$ starts to concentrate will depend on $\DFc$.
Since \(\UV(\DFc,\rrdc)\) does not change if we scale both \(\DFc\) and \(\rrdc\) by the same constant, we first fix the scale. A convenient way to do so is to assume \(\normG(\DFc)=1\). Recalling~\eqref{normGequiv}, we see that assuming \(\normG(\DFc)= 1\) amounts to starting with an arbitrary symmetric positive definite matrix \(\DFc^2\) and scaling it so that it is ``just under" \(\DFLG^2\). To present our theorem on the concentration of \(\PfL\), we need the following definition.
\begin{definition}[Bound on third order Taylor remainder]\label{def:del3-omega}
For a set \(\UV\) containing \(\map\) and a symmetric \(\DFL\succ0\), define
\begin{EQA}
\dltwb(\UV,\DFL) &=&\sup_{\thv\in\UV}\frac{|\lgd(\thv)-\lgd(\map)-\|\DFLG(\thv-\map)\|^2/2|}{\|\DFL(\thv-\map)\|^2/2}.\label{def:omega}
\end{EQA} 
\end{definition}
\begin{theorem}\label{thm:conc}Let \(\DFc\succ0\) satisfy \(\normG(\DFc)=1\). Suppose 
\begin{EQA}\label{n2lgd}
&&\nabla^2\lgd(\thv)\succeq\nabla^2\lgd(\map) - \DFc^2=\DFLG^2-\DFc^2,\qquad\forall \thv\in \R^{\dimp}.
\end{EQA} Furthermore, suppose \(\dltwb(\UV(\DFc,\rrdc),\DFc)\leq 1/3\), where \(\UV(\DFc,\rrdc)\) is as in~\eqref{def:UV-DFL}. Then 
\begin{EQA}
&&\PfL(\UV(\DFc,\rrdc)^c)\leq \frac13\exp\bigg(\!-\bigl(\rrdc-3\sqrt{\dimL(\DFc)}\bigr)^2/3\bigg),\qquad\forall\rrdc\geq3+3\sqrt{\dimL(\DFc)}.
\end{EQA} 
\end{theorem}
The proof of Theorem~\ref{thm:conc} is given in Appendix~\ref{app:conc}. Thus we see that the critical radius \(\rrdL\) such that \(\PfL\) concentrates in \(\UV(\DFc,\rrdc)\) is \(\rrdc\sim\sqrt{\dimL(\DFc)}\). This is precisely the same as for the LA, the Gaussian distribution \(\lap\) defined in~\eqref{df:lap}. Indeed, using well-known results on Gaussian concentration (see Corollary~\ref{corr:gaussconc}), we have
\begin{EQA}\label{gauss-tail-0}
\lap(\UV(\DFc,\rrdc)^c)&\leq &\exp\big(\!-\bigl(\rrdc-\sqrt{\dimL(\DFc)}\bigr)^2/2\big),\qquad\forall\rrdc\geq\sqrt{\dimL(\DFc)}.
\end{EQA} 
The takeaway message is that \(\PfL\) concentrates in the same elliptic sets that the Gaussian \(\lap\) concentrates in, though the radius \(\rrdL\) should be taken somewhat larger (by a factor of 3).  

\begin{remark}[Comparison to~\cite{spok23}]
Our result is similar to Theorem 4.1, (4.1) in~\cite{spok23}. However, there are two differences. First, we have written the tail bound in a way that makes explicit the decay rate with \(\rrdL\). Second, we have corrected two slight errors in the proof in that work, as explained in Remarks~\ref{rk:conc1} and~\ref{rk:conc2}.
\end{remark}

We now discuss why the conditions of Theorem~\ref{thm:conc} are needed, and when they are satisfied. 
The reason we need a Hessian lower bound is that \(\lgd\) should grow quadratically to obtain Gaussian-like tail decay of \(\PfL\propto \ex^{-\lgd}\). The form of the lower bound, \(\DFLG^2-\DFc^2\), is simply convenient for our calculations. Note that this is a relatively weak lower bound since \(\DFc^2\) is ``just under" \(\DFLG^2\) and therefore \(\DFLG^2-\DFc^2\) is only positive semi-definite. Indeed, since \(\normG(\DFc)=1\), the characterization~\eqref{normGequiv} gives \(\DFc^2\preceq\DFLG^2\) and \(\|\DFc\uv\|=\|\DFLG\uv\|\) for some \(\uv\neq0\). The specific way in which we use~\eqref{n2lgd} is through convexity of the function \(\lgd_0(\uv):=\lgd(\map+\uv)-\lgd(\map)-\frac12\uv^\T(\DFLG^2-\DFc^2)\uv\). Convexity implies\(\lgd_0\) grows linearly away from the boundary of \(\{\|\DFc\uv\|\leq\rrdc\}\). In turn, we can control \(\lgd_0\) on the boundary of this ball by using that \(\dltwb\leq1/3\), which implies \(\lgd_0\) grows quadratically away from zero in the local region.

To discuss when the conditions are satisfied, we consider the setting of Section~\ref{subsub:struct}, where \(\lgd=L+\|G\cdot\|^2/2\) and \(L\) is convex and, informally, of order \(\mathcal O(n)\). In this context, the condition~\eqref{n2lgd} can be rewritten as
\begin{EQA}\label{DFL-lb}
 \nabla^2L(\thv)&\succeq&\nabla^2L(\map)-\DFc^2\qquad\forall\thv\in\R^{\dimp}.
\end{EQA} Since \(\nabla^2L(\thv)\succeq0\) for all \(\thv\) this condition can be achieved by taking \(\DFc^2\succeq \nabla^2L(\map)\). Combining this with the fact that \(\normG(\DFc)=1\), and recalling~\eqref{normGequiv}, we can summarize a sufficient condition on \(\DFc\) as follows:
\begin{EQA}\label{DFL-mustbe}
\DFc^2&=&\nabla^2L(\map)+G_0^2,\qquad 0\preceq G_0^2\preceq G^2, \quad \exists\,\uv\;\text{s.t.} \;\|G_0\uv\|=\|G\uv\|.
\end{EQA} 
Note that when \(L\) is convex but not strongly convex,~\eqref{DFL-mustbe} is also \emph{necessary}. Indeed, we cannot choose \(\DFc^2\) any lower than \(\nabla^2L(\map)\) in the positive definite ordering, since there is no uniform lower bound on \(\nabla^2L\). 


Next, let us assume for ease of illustration that \(\nabla^2L(\map)\gtrsim  nI_{\dimp}\) and \(\DFc\) is as in~\eqref{DFL-mustbe}. Then \(\DFL\gtrsim\sqrt nI_{\dimp}\). (See Section~\ref{subsub:not} for the definition of \(A\gtrsim B\) for positive definite matrices \(A\) and \(B\).) To see when \(\dltwb\) is small, we first note the following inequality, shown in Appendix~\ref{app:subsec:taylor}:
\begin{EQA}\label{del3bds1}
&&\dltwb(\UV(\DFc,\rrdc),\DFc) \leq \frac{\rrdc}{3}\sup_{\|\DFc\uv\|\leq\rrdc}\|\nabla^3\lgd(\map+\uv)\|_{\DFc}.
\end{EQA} From the tensor norm definition in~\eqref{TDdef}, we see that \(\|\nabla^3\lgd\|_{\DFc}\) is like ``dividing" \(\nabla^3\lgd\) by three powers of \(\DFc\). Therefore, since \(\nabla^3\lgd=\nabla^3L\sim n\) and \(\DFc\) is at least of order \(n^{1/2}\), we get that \(\|\nabla^3\lgd\|_{\DFc}\sim n^{-1/2}\). Also, note that Theorem~\ref{thm:conc} requires that \(\rrdc\) is at least a constant multiple of \(\sqrt{\dimL(\DFc)}\). Thus for \(\rrdc\|\nabla^3\lgd\|_{\DFc}\) to be small (and therefore for \(\dltwb\) to be small, by~\eqref{del3bds1}), the quantity \(\sqrt{\dimL(\DFc)/n}\) should be sufficiently small. The condition \(\dimL(\DFc)\ll n\) is needed not only to ensure \(\UV(\DFc, \rrdc)\) contains most of the mass of \(\PfL\), but also to ensure this set is contained within a small ball around the MAP \(\map\). Indeed,  if \(\DFc\gtrsim \sqrt n I_{\dimp}\) and \(\rrdc=C\sqrt{\dimL(\DFc)}\), then the set \(\UV(\DFc, \rrdc)\) is contained within a ball of radius \(\mathcal O(\sqrt{\dimL(\DFc)/n})\) with respect to the standard Euclidean distance. Thus if \(\dimL(\DFc)\ll n\), then the posterior has the desirable property of concentrating in a shrinking region centered at the MAP \(\map\), reflecting low uncertainty in the parameter \(\thv\). 


Finally, let us return to the more realistic, anisotropic case mentioned at the end of Section~\ref{subsub:struct}, in which the smallest eigenvalues of \(\nabla^2L(\map)\) are much smaller than \(n\). In turn, the smallest eigenvalues of \(\DFc\) can be much smaller than \(\sqrt n\), since \(\DFc^2\) is allowed to be as small as \(\nabla^2L(\map)\). The effect of this is that the operator norm \(\|\nabla^3\lgd\|_{\DFc}\) may actually be much larger than \(n^{-1/2}\). For example, in Lemma~\ref{lma:compare} in the inverse problem setting, we show that \(\|\nabla^3\lgd\|_{\DFc}\) can be as large as \(\sqrt{\dimp/n}\). Similarly, the ellipse \(\UV(\DFc, \rrdc)\) need not be contained inside a small ball; it may be wide in certain directions.

\subsection{TV accuracy of the LA}\label{STV}The following lemma conveniently decomposes the TV distance into a ``local" TV distance over a set \(\UV\) and two tail probabilities. 
\begin{lemma}[Lemma 2.11 in~\cite{katsBVM}]\label{lma:TV}
For any Borel-measurable \(\UV\subseteq\R^{\dimp}\), we have
\begin{EQA}\label{TVdecomp}
&&\TV(\PfL,\lap)\leq \PfL(\UV^c) + \lap(\UV^c) + \TV(\PfL\vert_{\UV},\lap\vert_{\UV}),
\end{EQA} where \(\PfL\vert_{\UV},\lap\vert_{\UV}\) are the restrictions of \(\PfL,\lap\) to \(\UV\), respectively, normalized to be probability measures.
\end{lemma}
Note that here, \(\UV\) need not be of the form~\eqref{def:UV-DFL} studied in our concentration result. Lemma~\ref{lma:TV} already suggests that there is flexibility in the choice of \(\UV\). Our bound on \(\TV(\PfL\vert_{\UV},\lap\vert_{\UV})\) will introduce another degree of freedom, through a choice of matrix \(\DFL\). 
Recall the important Definition~\ref{def:effdim} of norm and effective dimension of a matrix \(\DFL\) relative to \(\DFLG\): \(\normG(\DFL)=\|\DFLG^{-1}\DFL\|\), and \(\dimL(\DFL)=\Tr(\DFLG^{-2}\DFL^2)/\normG(\DFL)^2\). We also introduce two more quantities.
\begin{EQA}
\dltwb_3(\UV,\DFL) &=&  \sup_{\thv\in\UV}\bigl\|\DFL^{-1}\bigl(\nabla^2\lgd(\thv)-\nabla^2\lgd(\map)\bigl)\DFL^{-1}\bigr\|,\label{def:omega3}\\
\dltwu_3(\UV,\DFL) &=& \sup_{\thv\in\UV}\bigl\|\DFL^{-1}\bigl(\nabla^2\lgd(\thv)-\nabla^2\lgd(\map)\bigl)\DFL^{-1}\bigr\|/\|\DFL(\thv-\map)\|.\label{def:tau3}
\end{EQA} The following is the strongest version of our result, which is based on a bound on \(\TV(\PfL\vert_{\UV},\lap\vert_{\UV})\), combined with Lemma~\ref{lma:TV}.
\begin{theorem}\label{app:thm:TV}
Let \(\UV\) be a set containing \(\map\), the minimizer of \(\lgd\). For any symmetric \(\DFL\succ0\) such that \(\lambda(\UV,\DFL) :=\normG(\DFL)^{-2}- \dltwb_3(\UV,\DFL)>0\), it holds
\begin{EQA}
\label{eq:TV}
	\TV(\PfL,\lap)\leq&& \frac{\dltwu_{3}(\UV,\DFL)}{4\sqrt{\lambda(\UV,\DFL)\lap(\UV)}}\left(\Tr(\DFLG^{-2}\DFL^2) + \|\DFLG^{-1}\DFL\|^2\right) + \PfL(\UV^c) + \lap(\UV^c).
\end{EQA}
\end{theorem}
See Section~\ref{sec:TVproof} for the proof, and Section~\ref{sec:comparegen} for a preview of the proof. Next, we make several simplifications. First, note that \(\dltwb_3(\UV,\DFL)\normG(\DFL)^2 = \dltwb_3(\UV, \DFL/\normG(\DFL))\) and \( \dltwu_3(\UV,\DFL)\normG(\DFL)^3 =  \dltwu_3(\UV, \DFL/\normG(\DFL))\), which is clear from the definitions~\eqref{def:omega3} and~\eqref{def:tau3}. It is then straightforward to show that both the condition \(\lambda(\UV,\DFL)>0\) and the righthand side~\eqref{eq:TV} can be expressed entirely in terms of the matrix \(\DFL/\normG(\DFL)\). Therefore we can assume without loss of generality that \(\normG(\DFL)=1\). The sum of the trace and norm in~\eqref{eq:TV} simplifies to \(\dimL(\DFL)+1\), which we bound by \(2\dimL(\DFL)\). Finally, we make the stronger assumption
\(\dltwb_3(\UV,\DFL)\leq 1/2\), which implies \(\lambda(\UV,\DFL)\geq1/2\), and we assume \(\lap(\UV)\geq1/2\). These simplifications lead to the following bound, which is coarser than~\eqref{eq:TV} only by a constant factor.
\begin{corollary}\label{thm:TV}
Let \(\UV\) be a set containing \(\map\) such that \(\lap(\UV)\geq1/2\). Then
\begin{EQA}\label{eq:TV:normalized}
\TV(\PfL,\lap)&\leq&\dltwu_3(\UV,\DFL)\dimL(\DFL)+\PfL(\UV^c) + \lap(\UV^c)\\
&\forall& \DFL\succ0:\;\normG(\DFL)=1,\;\dltwb_3(\UV,\DFL)\leq1/2.\label{cond:TV:normalized}
\end{EQA}
\end{corollary}
\begin{remark}
The only nontrivial condition needed for this result is that \(\dltwb_3(\UV,\DFL)\leq1/2\). In words, this means that in the set \(\UV\), the Hessian \(\nabla^2\lgd\) does not deviate too much from its value at \(\map\), where the deviation is measured in a \(\DFL\)-weighted operator norm.
\end{remark}
Before discussing our result, 
we make one further simplification. We show in Appendix~\ref{app:subsec:taylor} that
\begin{equation}\label{del3bds2}
\begin{split}
\dltwu_3(\UV,\DFL)&\leq\sup_{\thv\in\UV}\|\nabla^3\lgd(\thv)\|_{\DFL},\\
 \dltwb_3(\UV,\DFL)&\leq \sup_{\thv\in\UV}\|\nabla^3\lgd(\thv)\|_{\DFL}\sup_{\thv\in\UV}\|\DFL(\thv-\map)\|.
\end{split}\end{equation}
This gives the following most explicit (but coarsest) version of our bound, assuming again that \(\lap(\UV)\geq1/2\):
\begin{EQA}\label{TVsimpsimp}
\TV(\PfL,\lap)&\leq& \sup_{\thv\in\UV}\|\nabla^3\lgd(\thv)\|_{\DFL}\;\dimL(\DFL) + \PfL(\UV^c) + \lap(\UV^c),\\
 \forall \DFL&\succ&0:\;\normG(\DFL)=1,\quad\sup_{\thv\in\UV}\|\DFL(\thv-\map)\|\;\sup_{\thv\in\UV}\|\nabla^3\lgd(\thv)\|_{\DFL}\leq1/2.\label{TVsimpsimpcond}
\end{EQA} 
Thus the three factors influencing the bound are \(\lgd\), \(\DFL\), and \(\UV\), with the latter two being free parameters. We can think of \(\sup_{\thv\in\UV}\|\DFL(\thv-\map)\|\) as the ``radius" of the set \(\UV\) with respect to the distance defined by \(\DFL\). Although not strictly necessary, it is natural to define \(\UV\) to be a ball in the \(\DFL\) distance, i.e. \(\UV=\UV(\DFL, \rrdc)\), as in~\eqref{def:UV-DFL}. Then this radius is precisely \(\rrdc\) itself, and the second condition in~\eqref{TVsimpsimpcond} reduces to \(\rrdc\,\sup_{\|\DFL\uv\|\leq\rrdc}\|\nabla^3\lgd(\map+\uv)\|_{\DFL}\leq1/2\). 

With this choice of \(\UV\), our free parameters are now \(\DFL\) and \(\rrdc\). For any \(\DFL\), the optimal \(\rrdc\) is the smallest \(\rrdc\) such that \( \PfL(\UV(\DFL,\rrdc)^c) + \lap(\UV(\DFL,\rrdc)^c)\) is negligible compared to the first term on the righthand side of~\eqref{TVsimpsimp}. In particular, we should at least take \(\rrdc\geq\sqrt{\dimL(\DFL)}\). Indeed, by a standard result, for each \(\DFL\succ0\), the Gaussian distribution \(\lap\) concentrates in a thin shell around the sphere \(\{\thv\in\R^{\dimp}\,:\,\|\DFL(\thv-\map)\|=\sqrt{\dimL(\DFL)}\}\). How much larger \(\rrdc\) must be depends on the size of the first term in~\eqref{TVsimpsimp}, and on how rapidly the tails of \(\PfL\) decay. 
Once \(\rrdc=\rrdc(\DFL)\) has been chosen, one may then minimize the first part of the bound over \(\DFL\) satisfying \(\normG(\DFL)=1\) and \(\rrdc\,\sup_{\|\DFL\uv\|\leq\rrdc}\|\nabla^3\lgd(\map+\uv)\|_{\DFL}\leq1/2\).

For the sake of concreteness, we now state a full bound (including the tails) under the conditions of Theorem~\ref{thm:conc}, which guarantee \(\PfL\) has sub-Gaussian tail decay, and then discuss the optimal choice of \(\rrdc\) in this case. 
\begin{corollary}\label{corr:DD0} If $\mathrm{(1)}\,$\(\normG(\DFL)=1\), $\mathrm{(2)}\,$\(\nabla^2\lgd(\thv)\succeq\DFLG^2-\DFL^2\) for all \(\thv\in\R^{\dimp}\), and $\mathrm{(3)}\,$\(\rrdc\sup_{\|\DFL\uv\|\leq \rrdc}\|\nabla^3\lgd(\map+\uv)\|_{\DFL}\leq 1/2\) for some \(\rrdc\geq 3\sqrt{\dimL(\DFL)}+3\), then
\begin{EQA}\label{eq:TVsimpsimp}
\TV(\PfL,\lap) &\leq& \sup_{\|\DFL\uv\|\leq \rrdc}\|\nabla^3\lgd(\map+\uv)\|_{\DFL}\;\dimL(\DFL) + 2\ex^{-\bigl(\rrdL-3\sqrt{\dimL(\DFL)}\bigr)^2/3}.
\end{EQA}
\end{corollary}
Note that (1) and (3) ensures the conditions~\eqref{TVsimpsimpcond} are satisfied, and (3) also gives that \(\dltwb(\UV(\DFc,\rrdc),\DFc)\leq 1/3\). Thus (1),(2),(3) imply the conditions of Theorem~\ref{thm:conc}. We also used~\eqref{gauss-tail-0} to bound \(\lap(\UV^c)\), and in particular, it is clear that for the specified \(\rrdc\), we have \(\lap(\UV)\geq1/2\).

Let us discuss this result in the setting described in Section~\ref{subsub:struct}: \(\lgd=L+\|G\cdot\|^2/2\), where \(L\) is convex but not strongly convex and informally, \(L=\mathcal O(n)\) with respect to \(n\). Using the reasoning from Section~\ref{Sconc}, taking \(\DFL\) of the form~\eqref{DFL-mustbe} ensures (1) and (2 are satisfied, and we then have \(\|\nabla^3\lgd\|_{\DFL}\sim n^{-1/2}\). More specifically, consider the case in which \(\|\nabla^3\lgd\|_{\DFL}\geq n^{-1/2}\) (or a constant multiple of \(n^{-1/2}\)). Also, suppose \(\dimL(\DFL)\geq\log n\), which is satisfied for typical choices of \(\DFL\) in the inverse problem context; see Lemma~\ref{lma:compare}.  Then we can take \(\rrdc=C\sqrt{\dimL(\DFL)}\) for some large enough \(C\) to ensure that the second term in~\eqref{eq:TVsimpsimp} is negligible compared to the first term. Finally, consider condition (3) upon substituting \(\rrdc=C\sqrt{\dimL(\DFL)}\). We see that the upper bound~\eqref{eq:TVsimpsimp} being small is a stronger requirement than (3), so we can essentially neglect this condition.

Summarizing the above discussion, there are constants \(C_1, C_2\) such that for all \(\DFL^2\) of the form~\eqref{DFL-mustbe} such that \(\|\nabla^3\lgd\|_{\DFL}\geq n^{-1/2}\) and \(\dimL(\DFL)\geq\log n\), we have
\begin{EQA}\label{explicit}
&&\TV(\PfL,\lap) \lesssim \sup_{\|\DFL\uv\|\leq C_1\sqrt{\dimL(\DFL)}}\|\nabla^3\lgd(\map+\uv)\|_{\DFL}\,\dimL(\DFL),
\end{EQA} whenever this upper bound is smaller than \(C_2\). In particular, we at least need that \(n\gg\dimL(\DFL)^2\), since \(\|\nabla^3\lgd\|_{\DFL}\geq n^{-1/2}\). This is stronger than the condition \(n\gg\dimL(\DFL)\) needed for concentration in Section~\ref{Sconc}. 

The only free parameter remaining in~\eqref{explicit} is \(\DFL\). How to optimally choose \(\DFL\) depends on the specific structure of \(\lgd\), so we postpone the discussion of this choice to Section~\ref{inv} on the Bayesian inverse problem.


\subsection{Comparison to prior work}\label{sec:comparegen} We compare our results and assumptions to those of the three prior relevant works. We improve on Theorem 4.2 of~\cite{spok23} by decreasing the dependence on effective dimension from \(\dimL^{3/2}\) to \(\dimL\), and by introducing the optimization over \(\DFL\). Another important difference is that the TV bound presented in~\cite{spok23} relies on the concentration result requiring strong convexity of \(\lgd\). In contrast, strong convexity is \emph{not} required in Theorem~\ref{app:thm:TV} and Corollary~\ref{thm:TV}. We did use strong convexity in Corollary~\ref{corr:DD0}, but only to demonstrate \emph{one} way the tail probability \(\PfL(\UV^c)\) can be bounded.
The works~\cite{katsBVM} and~\cite{bp} have demonstrated that the tail probability can be made small without requiring \(\lgd\) to be strongly convex.

Next, we show that our results are strictly stronger than Theorem 2.2 of~\cite{katsBVM} and Theorem 17 of~\cite{bp}, and that these two theorems have the same form as our bounds. 
For simplicity, we neglect constants in the comparison, so \(\lesssim\) below denotes bounds that hold up to absolute constant.

~\cite{katsBVM} also starts from the decomposition~\eqref{TVdecomp}, and then further bounds both the local TV distance and the two tail probabilities.  As we show in Appendix~\ref{app:compare} the author's result can be written in the following form, under slightly stronger local conditions than in our Corollary~\ref{thm:TV}, together with tail conditions:
\begin{EQA}\label{katsBVM-simple-ub}
&&\TV(\PfL,\lap)\lesssim \dltwu_3(\UV,\DFLG)\;\dimL(\DFLG)  + \mathrm{UB}(\PfL(\UV^c)) + \mathrm{UB}(\lap(\UV^c)).
\end{EQA} Here, UB stands for upper bound, and in~\eqref{bp-simple-ub} below, LB stands for lower bound. The above bound is precisely as in our~\eqref{eq:TV:normalized}, except that the tails have been bounded, and the specific choice \(\DFL=\DFLG\) has been made.

Next, we compare our results to Theorem 17 in~\cite{bp}. We show in Appendix~\ref{app:compare} that the authors' result can be brought into the following form:
\begin{EQA}\label{bp-simple-ub}
\TV(\PfL,\lap)&\lesssim&\frac{\sup_{\thv\in\UV}\|\nabla^3\lgd(\thv)\|_{I_{\dimp}}}{\sqrt{\mathrm{LB}(\lambda(\UV,I_{\dimp}))\mathrm{LB}(\lap(\UV))}}\Tr(\DFLG^{-2}) +\mathrm{UB}(\PfL(\UV^c)) + \mathrm{UB}(\lap(\UV^c)).
\end{EQA}
This is precisely as in our Theorem~\ref{app:thm:TV} with the choice \(\DFL=I_{\dimp}\) and with the larger quantity \(\sup_{\thv\in\UV}\|\nabla^3\lgd(\thv)\|_{I_{\dimp}}\) instead of \(\dltwu_3(\UV,I_{\dimp})\), recalling~\eqref{del3bds2}. The bound~\eqref{bp-simple-ub} is shown to hold under six assumptions, one of which is a tail condition. Our single assumption that \(\lambda(\UV,\DFL) =\normG(\DFL)^{-2}- \dltwb_3(\UV,\DFL)>0\) is a slightly weaker version of the authors' Assumptions 3 and 5. We note that~\cite{bp} do not start with the decomposition from Lemma~\ref{lma:TV}, and some work is required to show that one of the tail terms in their TV bound is in fact an upper bound on \(\PfL(\UV^c)\). 

From~\eqref{bp-simple-ub}, one can do a similar normalization procedure as described below Theorem~\ref{app:thm:TV}; see Appendix~\ref{app:compare}. This leads to the following simplified form of the authors' bound, which differs from~\eqref{bp-simple-ub} only by a constant factor:
\begin{EQA}\label{bp-ub-norm}
\TV(\PfL,\lap)&\lesssim&\sup_{\thv\in\UV}\|\nabla^3\lgd(\thv)\|_{I_{\dimp}/\normG(I_{\dimp})}\;\dimL(I_{\dimp}/\normG(I_{\dimp})) + \mathrm{UB}(\PfL(\UV^c)) + \mathrm{UB}(\lap(\UV^c)).
\end{EQA}
This is precisely as in~\eqref{TVsimpsimp}, except that the tails have been bounded, and the specific choice \(\DFL=I_{\dimp}/\normG(I_{\dimp})\) has been made.

We highlight the two main reasons our results are stronger than those of~\cite{katsBVM} and~\cite{bp}. First, we allow for optimization over \(\DFL\). This leads to bounds that are tighter by an order of magnitude, as we will show in the inverse problem example. Alternatively, the flexibility in the choice of \(\DFL\) can be construed as weakening the assumptions: there need only exist \emph{some} \(\DFL\) for which~\eqref{cond:TV:normalized} is satisfied. Second, by leaving the tail probabilities as is, we have given the user the flexibility to bound them using any of the available techniques, and to choose which tail-related assumptions to impose. 

\paragraph{Proof comparison.}
Both in~\cite{bp} and~\cite{katsBVM}, as well as in the present work, the main lemma to bound the local contribution to the TV error is a version of the following: 
\begin{lemma}\label{lma:sobo}
Let \(\mu\propto \ex^{-f}\) and \(\nu\propto \ex^{-g}\) be two probability densities on a set \(\UV\subset \R^{\dimp}\). If \(\nabla^2f(\xv)\succeq\lambda\DFL^2\) for all \(\xv\in\UV\), then 
\begin{EQA}\label{TVLSI}
&&\TV(\mu,\nu)^2\leq \frac{1}{4\lambda}\E_{X\sim\nu}\left[\|\DFL^{-1}\nabla (f-g)(X)\|^2\right].
\end{EQA}
\end{lemma}
The work~\cite{bp} is the first to use this lemma, which is the key to improve the dimension dependence in the LA error bound from \(\sqrt{\dimp^3/n}\), obtained in prior works such as~\cite{helin2022non}, to \(\sqrt{\dimp^2/n}\). In~\cite{bp}, Lemma~\ref{lma:sobo} is stated and proved with \(\DFL=I_{\dimp}\). The lemma follows by combining Pinsker's inequality with the log Sobolev inequality. The latter applies because the Bakry-Emery criterion is satisfied, due to the uniform lower bound on the Hessian of \(f\) in \(\UV\). See Section D.3 of~\cite{bp} for references and more details. The work~\cite{katsBVM} then proves the lemma with \(\DFL=\DFLG\), simply by reducing it to the case \(\DFL=I_{\dimp}\) through a change of variables. But the proof in~\cite{katsBVM} (see that work's Lemma 2.11 and proof) goes through with exactly the same logic if we swap \(\DFLG\) for any other \(\DFL\). The key innovation in the present work is to use the lemma in its full generality, which allows us to take a step further and optimize over \emph{all possible \(\DFL\) satisfying the Hessian lower bound condition.}

\subsection{Proof of TV bound}\label{sec:TVproof}
%
%

In the below proof, it is helpful to note that \(\|\DFL^{-1}A\DFL^{-1}\|=\sup_{\|\DFL\uv\|=\|\DFL\vv\|=1}\langle A, \uv\otimes \vv\rangle\).
\begin{proof}[Proof of Theorem~\ref{app:thm:TV}]Assume without loss of generality that \(\map=0\). Due to Lemma~\ref{lma:TV}, it suffices to bound \(\TV(\PfL\vert_{\UV}, \lap\vert_{\UV})\). To do so, we apply Lemma~\ref{lma:sobo} with \(\mu=\PfL\vert_{\UV}\) and \(\nu=\lap\vert_{\UV}\propto\ex^{-g}\), with \(g(\xv)=\|\DFLG\xv\|^2/2\). First, we find \(\lambda\) such that \(\nabla^2\lgd\succeq\lambda\DFL^2\) uniformly over \(\UV\). Fix \(\xv\in\UV\). The definition~\eqref{def:omega3} of \(\dltwb_3=\dltwb_3(\UV,\DFL)\) implies 
\begin{EQA}
&&\nabla^2\lgd(\xv)\succeq \DFLG^2-\dltwb_3\DFL^2\succeq (\normG(\DFL)^{-2}-\dltwb_3)\DFL^2.\end{EQA}
Thus we can take \(\lambda=\normG(\DFL)^{-2}-\dltwb_{3}\). Next, to apply~\eqref{TVLSI}, we study \(\nabla (f-g)\). Note that \(\nabla g(\xv)=\nabla^2f(0)\xv\). Now, define \(h(t)=\langle\nabla f(t\xv),\uv\rangle\) for fixed \(\xv,\uv\). We then have
\begin{EQA}
\langle \nabla(f-g)(\xv),\uv\rangle &= &\langle\nabla f(\xv)-\nabla f(0) -\nabla^2f(0)\xv,\uv\rangle \\
&=& h(1)-h(0) - h'(0) = \int_0^1(h'(t)-h'(0))dt \\
&=& \int_0^1\langle\nabla^2f(t\xv)-\nabla^2f(0),\xv\otimes\uv\rangle dt.
\end{EQA} Thus
\begin{EQA}
\|\DFL^{-1}\nabla(f-g)(\xv)\|&=&\sup_{\|\DFL\uv\|=1}\langle \nabla(f-g)(\xv),\uv\rangle\\
&\leq &\sup_{\|\DFL\uv\|=1}\int_0^1\|\DFL^{-1}(\nabla^2f(t\xv)-\nabla^2f(0))\DFL^{-1}\|\|\DFL\xv\|\|\DFL\uv\|\\
&\leq &\dltwu_3\|\DFL\xv\|^2/2,\qquad\forall \xv\in\UV,\end{EQA} using the definition~\eqref{def:tau3} of \(\dltwu_3=\dltwu_3(\UV,\DFL)\).
We now return to~\eqref{TVLSI} with \(\nu=\lap\vert_{\UV}\) and note that for a function \(F\), it holds \(\E_{X\sim\lap\vert_{\UV}}[F(X)] = \E[F(\DFLG^{-1}Z)\Ind(\DFLG^{-1}Z\in\UV)]/\lap(\UV)\), where \(Z\sim\mathcal N(0, I_{\dimp})\). Using this observation and the above estimate gives

\begin{EQA}\label{Elapvert}
\E_{X\sim \lap\vert_{\UV}}&&\left[\left\|\DFL^{-1}\nabla(f-g)(X)\right\|^2\right]\leq \frac{\dltwu_{3}^2}{4\lap(\UV)}\E\left[\|\DFL\DFLG^{-1}Z\|^4\right]\\
&&\leq\frac{\dltwu_{3}^2}{4\lap(\UV)}\left(\Tr^2(\DFLG^{-1}\DFL^2\DFLG^{-1})+2\Tr((\DFLG^{-1}\DFL^2\DFLG^{-1})^2)\right).
\end{EQA}
In the second line, we used that \(\E[(\gaussv^{\T}B\gaussv)^2] = \Tr^2(B)+2\Tr(B^2)\) for a symmetric matrix \(B\), which we show below. Now, let \(B=\DFLG^{-1}\DFL^2\DFLG^{-1}\). Using that \(\Tr(B^2)\leq \|B\|\Tr(B)\) for a symmetric matrix \(B\succ0\), the expression in parentheses in the last line is \(\Tr^2(B)+2\Tr(B^2)\leq \Tr^2(B)+2\|B\|\Tr(B)\leq (\Tr(B)+\|B\|)^2\). Also, note that \(\|B\| = \|\DFLG^{-1}\DFL\|^2\). Applying Lemma~\ref{lma:sobo} now gives
\begin{EQA}
&& \TV(\PfL\vert_{\UV},\lap\vert_{\UV})^2\leq \frac{\dltwu_3^2}{16\lambda\lap(\UV)}\left(\Tr(\DFLG^{-2}\DFL^2)+\|\DFLG^{-1}\DFL\|^2\right)^2.
\end{EQA}
Taking the square root and using Lemma~\ref{lma:TV} concludes the proof.

To show \(\E[(\gaussv^{\T}B\gaussv)^2] = \Tr^2(B)+2\Tr(B^2)\), assume without loss of generality that \(B\) is diagonal with diagonal entries \(\eta_j\). Then
\begin{EQA}
\E[(\gaussv^{\T}B\gaussv)^2] &=&\E[(\sum_i\eta_i\gaussv_i^2)^2] =\sum_{i\neq j}\eta_i\eta_j +3\sum_i\eta_i^2 \\
&=& (\sum_i\eta_i)^2 + 2\sum_i\eta_i^2=\Tr^2(B)+2\Tr(B^2).
\end{EQA}

\end{proof}

%
%
%

\section{Application to generalized linear inverse problem}\label{inv}

%
%

We describe the set-up in Section~\ref{sec:inv:setup}, including assumptions on the forward operator and the statistical model. In Section~\ref{sub:opt:inv} we present our TV bound using the optimized choice of \(\DFL\), and compare to bounds from prior works, which do not choose \(\DFL\) optimally. We give more details about the optimization over \(\DFL\) and computation of the TV bound in Section~\ref{subsec:gam0}. All omitted proofs are in Appendix~\ref{app:inv}, unless indicated otherwise.
\subsection{Set-up}\label{sec:inv:setup}
Briefly, the set-up is that we observe noisy perturbations of \((\calR q^*)(j/n)\), \(j=1,\dots,n\), and we wish to recover the first \(\dimp\) coefficients of an orthonormal basis representation of \(q^*\). The operator \(\calR\) and associated basis functions are introduced in Section~\ref{subsub:R}, and the statistical model is presented in Section~\ref{subsub:stat}. First, we introduce some notation. 

\subsubsection{Notation}\label{subsub:not} Suppose \(a,b>0\). Then  \(a\lesssim b\) and \(a\asymp b\) means there exist constants \(0<c_1<c_2\) such that \(a\leq c_2b\) and \(c_1b\leq a\leq c_2b\), respectively. Similarly, for two symmetric positive definite matrices \(A,B\), we say \(A\lesssim B\) and \(A\asymp B\) if there exist constants \(0<c_1<c_2\) such that \(\uv^\T A\uv \leq c_2\uv^\T B\uv\) and \(c_1\uv^\T B\uv\leq \uv^\T A\uv \leq c_2\uv^\T B\uv\) for all \(\uv\in\R^{\dimp}\), respectively. Finally, for two symmetric matrices \(A,B\) such that \(B\succ0\), we write \(|A|\lesssim B\) to denote \(|\uv^\T A\uv|\leq c_2\uv^\T B\uv\) for all \(\uv\in\R^{\dimp}\) and some constant \(c_2>0\).

The suppressed constants \(c_1,c_2\) can depend on all problem parameters except \(\dimp\), \(n\), and any quantities derived from \(\dimp\) and \(n\).


\subsubsection{Linear operator and basis functions}\label{subsub:R}
Let \(\calR:E\subseteq L^2[0,1]\to L^2[0,1]\) be a compact linear operator, where \(E\) is a linear subspace of \(L^2[0,1]\). We let \(\phi_k\), \(k=1,2,3,\dots\) be a complete orthonormal basis of \(E\) consisting of eigenfunctions of \(\calR^\T\calR\). The corresponding eigenvalues are denoted \(\lambda_k\), \(k=1,2,3,\dots\). We assume \(\lambda_1\geq\lambda_2\geq\dots>0\), and that the eigenvalues decrease to zero as
\begin{EQA}\label{lambda2beta}
&&\lambda_k\asymp k^{-2\beta},\qquad k=1,2,3,\dots,\qquad\beta\geq0.
\end{EQA}
This means there exist constants \(0<c<C\) such that \(ck^{-2\beta}\leq\lambda_k\leq Ck^{-2\beta}\) for all \(k=1,2,3,\dots\). Note that the functions
\begin{EQA}\label{phi-psi}
\psi_k&=&\calR\phi_k/\sqrt{\lambda_k},\qquad k=1,2,3,\dots
\end{EQA}
form a complete orthonormal basis of the range space of \(\calR\); they are also the eigenfunctions of \(\calR\calR^\T\). The main point is that \(\calR\) is a \emph{smoothing} operator, sending \(\phi_k\) to \(k^{-\beta}\psi_k\). Therefore, if e.g. the \(\phi_k\) and \(\psi_k\) are both an orthonormal basis of sinusoids, then \(\calR\) has the effect of damping high frequencies.

\begin{assumption}[Boundedness in \(L^\infty\) norm]\label{assume:bdd}
We assume that \(\sup_{k=1,2,3,\dots}\|\psi_k\|_\infty \lesssim 1\).
\end{assumption}
This assumption is used in Lemma~\ref{lma:effdim:radon} to bound a third derivative operator norm. We also use it as an auxiliary tool to verify the next assumption.
\begin{assumption}[Near orthogonal discretized basis]\label{assume:ortho}Let \(\Psi\in\R^{\dimp\times\dimp}\) with entries \(\Psi_{k\ell}=\sum_{j=1}^n\psi_k(j/n)\psi_\ell(j/n)\). There exists \(\lambda\in\R\) independent of \(\dimp\) and \(n\) such that
\begin{EQA}\label{Ckbg}
&&|\Psi - nI_{\dimp}|\precsim \diag(k^{\lambda}),
\end{EQA} where the suppressed constant may depend on \(\lambda\).\end{assumption} In words, we assume that discretizing the first \(\dimp\) basis functions \(\psi_k\), \(k=1,\dots,\dimp\) to vectors \((\psi_k(j/n))_{j=1}^n\) in \(\R^n\) nearly preserves orthonormality. See the beginning of Section~\ref{sub:opt:inv} for an explanation of how this assumption is used. Note that~\eqref{Ckbg} is reasonable because \(n^{-1}\Psi_{k,\ell}\) is the Riemann sum approximation of \(\int_0^1\psi_k(x)\psi_\ell(x)dx=\delta_{k,\ell}\). Choosing \(\lambda>0\) in the bounds~\eqref{Ckbg} allows the Riemann error to grow as \(k,\ell\) become large. This room for error is useful when the derivatives of \(\psi_k\) grow with index. The following lemma gives sufficient conditions for the assumption to be satisfied. 
\begin{lemma}\label{lma:D:prelim}
If \(\sup_k\|\psi_k\|_\infty\lesssim 1\) (i.e. Assumption~\ref{assume:bdd} is satisfied) and if there is \(\delta\geq0\) such that the \(\psi_k\) satisfy \(\sup_{k=1,2,3,\dots}k^{-\delta}\|\psi_k'\|_\infty\les1\), then Assumption~\ref{assume:ortho} is satisfied for all \(\lambda>1+2\delta\).\end{lemma}
We now state a key result showing that the conditions of the lemma are verified for a large class of operators \(\calR\), with \(\delta=1\).
 \begin{proposition}\label{thm:psik0}
Let \(a\in C^2[0,1]\), \(b\in C^1[0,1]\), and assume there is \(a_0>0\) such that \(a(x)>a_0>0\) for all \(x\in[0,1]\). Let \(\calR\) be the operator on \(L^2[0,1]\) defined by
\begin{EQA}\label{calR-def0}
&&g=\calR f \quad \Leftrightarrow\quad ag' + bg = f,\quad g(0)=0,
\end{EQA} which is equivalent to the operator
\begin{EQA}\label{calR-int-def0}
&&\calR f(x)= e^{-C(x)}\int_0^x \frac{e^{C(t)}}{a(t)}f(t)dt,
\end{EQA}
where \(C\) is an antiderivative of \(b/a\). Then \(\calR\) is compact. Let \(\psi_k\), \(k=1,2,3,\dots\) be a complete orthonormal basis of \(L^2[0,1]\) consisting of eigenfunctions of \(\calR\calR^\T\), and let \(\lambda_k\) be the corresponding eigenvalues. Then \(\lambda_k\asymp k^{-2}\), and
\begin{EQA}\label{psik-assume0}
&&\|\psi_k\|_\infty \leq \Psi_0,\qquad \|\psi_k'\|_\infty \leq \Psi_1 k,
\end{EQA}
where the constants \(\Psi_0,\Psi_1\) depend only on the functions \(a\) and \(b\). Therefore, Assumptions~\ref{assume:bdd} and~\ref{assume:ortho} are both satisfied, for any \(\lambda>1+2\delta=3\).
\end{proposition}
See Appendix~\ref{app:eigen} for the proof. 
\begin{remark}
Recall from the beginning of the section that we started with a complete orthonormal basis \((\phi_k)_k\) of \(L^2[0,1]\) consisting of eigenfunctions of \(\calR^\T\calR\), and defined the \(\psi_k\) as \(\psi_k=\calR\phi_k/\sqrt{\lambda_k}\). In contrast, in Proposition~\ref{thm:psik0}, we defined \((\psi_k)_k\) as an arbitrary complete orthonormal basis of \(L^2[0,1]\) of eigenfunctions of \(\calR\calR^\T\). This is fine because the nullspaces of \(\calR\) and \(\calR^\T\) are trivial, which implies any basis \((\psi_k)_k\) of the latter type can be written as \(\psi_k=\calR\phi_k/\sqrt{\lambda_k}\), for a basis \((\phi_k)_k\) of the former type.  
\end{remark}

Taking \(a\equiv1\) and \(b\equiv0\) gives the Volterra operator \((\calR f)(x) = \int_0^xf(t)dt\), which was studied in a similar but linear Bayesian inverse problem context in~\cite{knapik2011bayesian}.


\begin{remark}[Equivalent formulation] It is straightforward to show that~\eqref{calR-int-def0} is equivalent to 
\begin{EQA}
&&\calR f(x)= S(x)\int_0^x K(t)f(t)dt,
\end{EQA}
for sufficiently smooth functions \(S,K\) of constant sign. Indeed, taking \(a=1/(SK)\) and \(b=-S'/(S^2K)\) gives back~\eqref{calR-int-def0}.
\end{remark}
\begin{remark}[Proof technique]We show in Appendix~\ref{app:eigen} that for the above \(\calR\), the operator \(\calR\calR^\T\) is the inverse of a Sturm-Liouville differential operator; see Lemma~\ref{lma:RRT}. Thus we can study the eigenvalues \(\lambda_k\) and eigenfunctions \(\psi_k\) of  \(\calR\calR^\T\) using the powerful theory of eigenvalue and eigenfunction asymptotics for Sturm-Liouville problems. See Appendix~\ref{app:eigen} for references and the proof. 
\end{remark}
\begin{remark}[Beyond $\calR$ of the type~\eqref{calR-int-def0}]\label{rk:basis}
Let \(\calS\) be any other compact operator such that \(\calS\calS^\T\) commutes with \(\calR\calR^\T\). This is equivalent to the inverse of \(\calS\calS^\T\) commuting with the Sturm-Liouville differential operator from Lemma~\ref{lma:RRT}. Then by standard theory, \(\calS\calS^\T\) shares an orthonormal basis of eigenfunctions \(\psi_k\) with \(\calR\calR^\T\); see e.g. Theorem 7, Chapter 28 of~\cite{lax2014functional}. Thus~\eqref{psik-assume0} is automatically satisfied. Note however that \(\calS\calS^\T\) need not have the same eigenvalues as \(\calR\calR^\T\).

\end{remark}

\subsubsection{Statistical model}\label{subsub:stat}
Define a one-parameter exponential family \(\rho(\cdot\mid\cdot)\) by
\begin{EQA}\label{exp-fam-def}
\rho(dy\mid s) &=& \e\l(sy - h(s)\r)\mu(dy),\quad y\in\mathcal Y\subset\R
\end{EQA} for some base measure \(\mu\) on \(\mathcal Y\). Note that \(h\) is the log normalization constant:
\begin{EQA}
h(s)&=&\log\int_{\mathcal Y} e^{sy}\mu(dy).
\end{EQA} As described in Section~\ref{intro:inv}, our motivating example is an inverse problem with Poisson noise. Thus \(\rho\) is a Poisson distribution, for which \(\mathcal Y=\{0,1,2,\dots\}\). The pmf of a Poisson$(\lambda)$ can be brought into the form~\eqref{exp-fam-def}, with \(h(s)=e^s=\lambda\). 

We assume standard regularity conditions on the exponential family ensuring that \(h\) is strictly convex and infinitely differentiable in the domain where it is finite, and furthermore, we assume this domain is all of \(\R\). See Chapters 7, 8 of~\cite{barndorff2014information}. The Poisson, binomial, and geometric distributions mentioned in Section~\ref{intro:inv} all satisfy these properties. 

Next, let $q^*\in L^2[0,1]$ be the ground truth. 
We assume that we are given independent observations $Y_1,\dots, Y_n$ from the distribution
\begin{EQA}\label{observations}
Y_j &\sim &\rho(\cdot \mid (\calR q^*)(j/n)),\quad j=1,\dots,n.\end{EQA} 
\begin{remark}By standard properties of exponential families, the mean of \(\rho(\cdot\mid s)\) is \(h'(s)\). Thus \(\E[Y_j] = h'((\calR q^*)(j/n))\). This means we can consider the \(Y_j\) to be noisy perturbations of \(h'((\calR q^*)(j/n))\), though the ``perturbation" is nonlinear. Since \(h''>0\), the function \(h'\) is monotone increasing and hence invertible. Thus observations of \(Y_j\) give information about \((\calR q^*)(j/n)\). We note that the model~\eqref{observations} can be considered a nonlinear analogue of the standard regression model with additive noise: \(Y_j = (\calR q^*)(j/n)+\xi_j\), \(j=1,\dots,n,\) for i.i.d. standard Gaussian \(\xi_j\).
\end{remark}
Next, we introduce the finite-dimensional basis approximation to \(q^*\).
\begin{definition}\label{def:qth-rhoth}For a given \(\thv\in\R^{\dimp}\), let \(q_{\thv}=\sum_{k=1}^{\dimp}\theta_k\phi_k\).
\end{definition}
We approximate \(q^*\) by \(q_{\thv}\) for an unknown \(\thv\in\R^{\dimp}\), and do inference over \(\thv\). In other words, the likelihood which we now specify is based on the approximate model
\begin{EQA}\label{model}
Y_j &\sim &\rho(\cdot \mid (\calR q_{\thv})(j/n)),\quad \thv\in \R^{\dimp},\quad j=1,\dots,n,
\end{EQA} even though the true data-generating process is~\eqref{observations}.
Note that we may write \( (\calR q_{\thv})(j/n)=R_j^\T\thv \), where \(R_1^\T,\dots,R_{n}^\T\in\R^{\dimp}\) are the rows of the \(n\times\dimp\) matrix \(R\) with entries
\begin{EQA}\label{Rjrows}
R_{jk}&=& (\calR \phi_k)(j/n),\quad k=1,\dots,\dimp,\quad j=1,\dots,n.
\end{EQA}
Using the approximate model~\eqref{model} and the definition~\eqref{exp-fam-def} of \(\rho\), we obtain the following negative log likelihood. We also record its Hessian.
\begin{EQA}\label{LL-radon}
L(\thv) &=& \sum_{j=1}^n\bigl[h(R_j^\T\thv) -Y_jR_j^\T\thv\bigr].\\
\nabla^2L(\thv)&=&\sum_{j=1}^nh''(R_j^\T\thv)R_jR_j^\T.\label{LL-hess-radon}
\end{EQA} 
We now consider the posterior of \(\thv\), given a Gaussian prior \(\thv\sim\ND(0,G^{-2})\), where 
\begin{EQA}\label{G2gamma}
G^2&=&\diag(k^{2\gamma})_{k=1}^{\dimp},\qquad\gamma\geq 0.
\end{EQA}
The negative log posterior is thus
\begin{EQA}\label{lgd:inv}
\lgd(\thv)&=&L(\thv) + \frac{1}{2}\|G\thv\|^2.
\end{EQA} 
Note that \(L\) is convex because \(h\) is convex. However, \(L\) is typically not strongly convex, because \(h\) is not. For example, for our motivating example of a Poisson distribution, \(h(s)=e^s\) is not strongly convex on \(\R\). Thus this setup is exactly as described in Section~\ref{subsub:struct}.

Since \(\lgd\) is strongly convex, it has a unique global minimizer \(\map :=\argmin_{\thv\in\R^{\dimp}}\lgd(\thv)\). 
The following is a key assumption about \(\map\).
\begin{assumption}\label{assume:Rqinfty}
There is a constant \(C\) independent of \(n\) and \(\dimp\) so that \(\|\calR q_{\map}\|_\infty\leq C\).
\end{assumption} We explain the reason for this assumption below Lemma~\ref{lma:hess}. Since \(\map\) is random (due to the randomness of the $Y_j$'s), the assumption should be interpreted to mean that we work on the event on which \(\|\calR q_{\map}\|_\infty\leq C\) is satisfied. For this assumption to be meaningful, the problem parameters should be such that this event has high probability. Now, since \(\calR\) is a smoothing operator, it is very reasonable for \(\|\calR q^*\|_\infty\) to be bounded. Therefore, boundedness of \(\|\calR q_{\map}\|_\infty\) follows if \(\calR q_{\map}\) is not too far from \(\calR q^*\). Thus the natural way to verify Assumption~\ref{assume:Rqinfty} is to solve the frequentist estimation problem of showing \(\calR q_{\map}\) is close to \(\calR q^*\) with high probability. This is a separate problem which lies outside the scope of the present work, which is why we consider Assumption~\ref{assume:Rqinfty} as given, rather than proving it. To give an example of how the assumption constrains the problem parameters, we note that if the regularization is constant (i.e. \(\gamma=0\)), then \(\dimp\) likely cannot be arbitrarily large relative to \(n\). On the other hand, if \(\gamma\) is well-adapted to the smoothness of \(q^*\), then classical theory of nonparametric estimation in linear settings suggests \(\dimp\) can indeed be arbitrarily large. See~\cite{knapik2011bayesian,cavalier2008nonparametric,tsybakov2010nonparametric}.

It is clear that \(\|\calR q_{\map}\|_\infty\leq C\) is a minimal requirement for inference to be meaningful, so in this work, we simply assume the user has appropriately specified the prior and the basis discretization level \(\dimp\) to achieve this. 

\subsection{Optimized TV bound and comparison to prior work}\label{sub:opt:inv}
We now present our  TV bound using an ``optimized" choice of \(\DFL\). Consider the following family of matrices \(\DFL^2\):
\begin{EQA}\label{DFLgam0def}
&&\DFL^2=\DFL^2(\gamma_0):=\nabla^2L(\map) + \diag(k^{2\gamma_0}),\qquad\gamma_0 \leq \gamma. 
\end{EQA} Note that this family includes \(\DFLG^2\) itself: \(\DFL^2(\gamma)=\DFLG^2\). For each fixed \(\gamma_0\), we apply Corollary~\ref{corr:DD0} with \(\DFL=\DFL(\gamma_0)\), and choose a radius \(\rrdc\) which makes the second term in the bound negligible compared to the first. To apply the corollary, we must compute \(\|\nabla^3\lgd\|_{\DFL(\gamma_0)}\) and \(\dimL(\DFL(\gamma_0))\). This generally challenging calculation is significantly simplified by the fact that \(\DFL^2(\gamma_0)\asymp\diag(n/k^{2\beta}+k^{2\gamma_0})\). We show this in Section~\ref{subsec:gam0} using the crucial Assumption~\ref{assume:ortho}. This is the main reason for the assumption.
\begin{theorem}\label{corr:TV:radon}Suppose Assumptions~\ref{assume:bdd},~\ref{assume:ortho},~\ref{assume:Rqinfty} are satisfied, and recall \(\lambda\) from Assumption~\ref{assume:ortho}. Define
\begin{EQA}\label{SdimL-def}
S_{\dimL}(\gamma_0)&:=&\sum_{k=1}^{\dimp}\frac{n+k^{2\gamma_0+2\beta}}{n+k^{2\gamma+2\beta}},\\
S_{\dltwu}(\gamma_0)&:= & \left( \sum_{k=1}^{\dimp}\frac{1}{n+k^{2\gamma_0+2\beta}}\right)^{1/2}.\label{Stau-def}
\end{EQA}
Then for all \(\gamma_0\leq\gamma\) such that \(2\gamma_0+2\beta\geq\lambda\) and \(S_{\dltwu}(\gamma_0)S_{\dimL}(\gamma_0)\leq c\) for a sufficiently small constant \(c\) independent of \(\gamma_0\), it holds
\begin{EQA}\label{TV-radon-2}
&&\TV( \PfL,\,\lap)\les S_{\dltwu}(\gamma_0)S_{\dimL}(\gamma_0).
\end{EQA}
\end{theorem}
The bound~\eqref{TV-radon-2} stems from~\eqref{eq:TVsimpsimp} and the fact that \(\|\nabla^3\lgd\|_{\DFL(\gamma_0)}\les S_{\dltwu}(\gamma_0)\) and \(\dimL(\DFL(\gamma_0))\asymp S_{\dimL}(\gamma_0)\), which we show in Section~\ref{subsec:gam0}. Note that \(1\leq S_{\dimL}(\gamma_0)\leq\dimp\) and \(S_{\dltwu}(\gamma_0)\leq\sqrt{\dimp/n}\) for all \(\gamma_0\), so the bound~\eqref{TV-radon-2} is guaranteed to be small if \(n\gg\dimp^{3/2}\). However, by optimizing the product \(S_{\dltwu}(\gamma_0)S_{\dimL}(\gamma_0)\) over \(\gamma_0\), we can significantly improve this bound. It is clear from the definitions~\eqref{SdimL-def} and~\eqref{Stau-def} that increasing \(\gamma_0\) increases \(S_{\dimL}\) but decreases \(S_{\dltwu}\). We show in Section~\ref{subsec:gam0} that the optimal choice of \(\gamma_0\) is just under \(\gamma-1/2\). This leads to the following, main result of Section~\ref{inv}.


\begin{corollary}\label{corr:TVgam0star}Suppose Assumptions~\ref{assume:bdd},~\ref{assume:ortho},~\ref{assume:Rqinfty} are satisfied and let \(\lambda\) be as in Assumption~\ref{assume:ortho}. Suppose \(2\gamma+2\beta>\max(2,\lambda+1)\), and let \(m=n^{1/(2\beta+2\gamma)}\), and \(m_0^*=n^{1/(2\beta+2\gamma_0^*)}\), where \(\gamma_0^*=\gamma-1/2-1/2m\). Then for all \(n\) large enough, 
\begin{EQA}
\TV(\PfL,\lap) &\les& \frac{(m\wedge\dimp)\sqrt{m_0^*\wedge\dimp}}{\sqrt n}\label{TVgam0} 
\end{EQA}whenever the upper bound is sufficiently small.
\end{corollary}
\begin{remark}When \(\dimp\geq n^{1/(2\beta+2\gamma_0^*)}\), the bound reduces to $$\TV(\PfL,\lap)\les \sqrt{n^{2/(2\beta+2\gamma)+1/(2\beta+2\gamma_0^*)-1}}.$$ Thus the bound is independent of \(\dimp\) for \(\dimp\) large enough. Furthermore, using the approximation \(\gamma_0^*\approx\gamma-1/2\), we get that the power of \(n\) is negative when \(2\beta+2\gamma\gtrapprox 2+\sqrt2\).
\end{remark}

\begin{remark}Proposition~\ref{thm:psik0} shows that for the family~\eqref{calR-int-def0} of $\calR$'s, the eigenvalues of \(\calR^\T\calR\) decay as \(k^{-2\beta}\) with \(\beta=1\), and that Assumption~\ref{assume:ortho} is satisfied for any \(\lambda>3\). Thus the condition on \(2\gamma+2\beta\) in Corollary~\ref{corr:TVgam0star} is that \(2\gamma+2=2\gamma+2\beta>4\). Therefore, the regularization parameter \(\gamma\) should satisfy \(\gamma>1\).
\end{remark}
We now show that our full TV bound~\eqref{TVgam0} is smaller by an order of magnitude than the \emph{local contribution} to the TV bounds of~\cite{bp} and~\cite{katsBVM}, and therefore smaller than their full TV bounds (consisting of the sum of the local and tail contributions). Recall from the discussion in Section~\ref{sec:comparegen} that in all three works, the local contribution to the TV bound is given by 
$$
\dimL(\DFL)\sup_{\thv\in\UV}\|\nabla^3\lgd(\thv)\|_{\DFL},
$$
for different choices of \(\DFL\). To be more precise, two of the three works use the Lipschitz constant of \(\nabla^2\lgd\) in place of the larger quantity \(\|\nabla^3\lgd(\thv)\|_{\DFL}\). Here, we work with the operator norm for simplicity. 

In order to facilitate the comparison, we assume for convenience that \(h'''\) is uniformly bounded over \(\R\). This allows us to extend the supremum of the derivative operator norm from \(\UV\) to all of \(\R^{\dimp}\). As a result, we can focus on the effect of different choices of \(\DFL\) on the local TV bounds, and neglect differences due to the set \(\UV\). Note that the above Theorem~\ref{corr:TV:radon} and Corollary~\ref{corr:TVgam0star} do \emph{not} require \(h'''\) to be uniformly bounded.

In the following lemma, we use the notation
\begin{EQA}
&&\|\nabla^3\lgd\|_{\DFL,\infty} :=\sup_{\thv\in\R^{\dimp}}\|\nabla^3\lgd(\thv)\|_{\DFL}.
\end{EQA}

\begin{lemma}\label{lma:compare}Suppose Assumptions~\ref{assume:bdd},~\ref{assume:ortho},~\ref{assume:Rqinfty} are satisfied and \(\|h'''\|_\infty<\infty\).  Let \(m,m_0^*,\gamma_0^*\) be defined as in Corollary~\ref{corr:TVgam0star}. Suppose \(\beta>1/2\), \(\gamma>1/2\), and \(2\beta+2\gamma>\max(2,\lambda+1)\). Then for all \(n\) large enough,
\begin{EQA}[rclrcl]
\dimL(\DFL(\gamma_0^*))&\asymp & m\wedge\dimp,\qquad &\|\nabla^3\lgd\|_{\DFL(\gamma_0^*),\infty}&\les& \sqrt{\frac{m_0^*\wedge\dimp}{n}},\\
\dimL(\DFLG)&=& \dimp,\qquad &\|\nabla^3\lgd\|_{\DFLG,\infty} &\les& \sqrt{ \frac{m\wedge\dimp}{n}},\label{tau3D2}\\
 \dimL(I_{\dimp}/\normG(I_{\dimp}))&\asymp& m\wedge\dimp,\qquad &\|\nabla^3\lgd\|_{I_{\dimp}/\normG(I_{\dimp}),\infty} &\les& \sqrt{\frac{(m\wedge\dimp)^{6\beta}}{n}}.
\end{EQA}
\end{lemma}
\begin{remark}[Tightness of operator norm bounds]
 Note that we have computed \(\dimL(\DFL)\) exactly, up to constant factors, but we have only obtained upper bounds on \(\|\nabla^3\lgd\|_{\DFL,\infty}\). To fairly compare the three different choices of \(\DFL\), we need to be sure the degree to which we have coarsened the operator norms is the same for all three \(\DFL\)'s. We confirm this in Appendix~\ref{app:tight}, under the simplifying assumption that the \(\psi_k\) are a basis of cosine functions. We show that after a first step in which we pass an absolute value inside a sum and bound all factors \(|h'''(t)|\) by \(\|h'''\|_\infty\), the rest of the bound is tight up to constant factors, for each of the three operator norms. 
\end{remark}

It is straightforward to check that in the regime \(2\beta+2\gamma>\max(2,\lambda+1)\) and \(\beta>1/2\), we have \(m_0^*\leq m^{6\beta}\) as soon as \(n\geq3^{2\beta+2\gamma}\). Thus the three bounds on \(\|\nabla^3\lgd\|_{\DFL,\infty}\) are ordered as follows:
\begin{EQA}
&& \sqrt{ \frac{m\wedge\dimp}{n}} \leq  \sqrt{\frac{m_0^*\wedge\dimp}{n}} \leq  \sqrt{\frac{(m\wedge\dimp)^{6\beta}}{n}}.
\end{EQA}
From here, the comparison between \(\DFL=\DFL(\gamma_0^*)\) and \(\DFL=I_{\dimp}/\normG(I_{\dimp})\) is clear: they both have the same \(\dimL(\DFL)\), but the latter has a larger value of \(\|\nabla^3\lgd\|_{\DFL,\infty}\), and therefore the latter leads to a larger product. Comparing \(\DFL=\DFL(\gamma_0^*)\) to \(\DFL=\DFLG\), we see that the latter has a larger \(\dimL(\DFL)\) but a smaller \(\|\nabla^3\lgd\|_{\DFL,\infty}\). However, the product for \(\DFL(\gamma_0^*)\) is smaller than the product for \(\DFLG\).

%

Table~\ref{table:cases} shows the squared TV bound \(\dimL(\DFL)^2\mathrm{UB}(\|\nabla^3\lgd\|_{\DFL,\infty})^2\), for the three choices of \(\DFL\) and for various values of \(\dimp\). Here, \(\mathrm{UB}(\|\nabla^3\lgd\|_{\DFL,\infty})\) denotes the operator norm upper bounds in Lemma~\ref{lma:compare}. Note that for our optimized choice, \(\DFL=\DFL(\gamma_0^*)\), the bound is dimension free: \emph{it does not depend on \(\dimp\) for all \(\dimp\) large enough.} Therefore, \(n\) can be arbitrarily large relative to \(\dimp\). The same is true for the choice \(\DFL=I_{\dimp}/\normG(I_{\dimp})\), though the bound is still larger than using the optimized \(\DFL(\gamma_0^*)\). The comparison between the bounds is shown in the 
last two columns in the table. 
In the second-to-last column, the improvement of \(\DFL(\gamma_0^*)\) over \(\DFLG\) is clear: the ratio is always larger than 1. In the last column, since \(\beta>1/2\), the improvement in the \(\dimp<m\) case is at least \(\dimp^2\). The improvement in the \(m<\dimp<m_0^*\) and \(m_0^*<\dimp\) cases are both at least \(m^{6\beta}/m_0^*\). Furthermore, assuming \(1/m\) is negligible and using that \(\beta+\gamma\geq1\), we have
\begin{EQA}
&&m_0^*=m^{\frac{\beta+\gamma}{\beta+\gamma_0^*}}\approx m^{\frac{\beta+\gamma}{\beta+\gamma-1/2}} \leq m^{2}
\end{EQA} Thus the improvement in the cases \(m<\dimp<m_0^*\) and \(\dimp>m_0^*\) is at least \(m^{6\beta-2}\).

Table~\ref{table:cases} justifies our statement in the introduction: our optimized \(\DFL\) gives a bound which is an order of magnitude tighter than the known bounds.

\begin{remark}[The bound of~\cite{katskew}]\label{rk:katskew}
Recall from the introduction that~\cite{katskew} derives a leading order decomposition of \(\TV(\PfL,\lap)\). This gives rise to an alternative TV upper bound, of a different form than the bounds we have studied so far. The upper bound requires a term called \(\epsilon_4(s)^2\) to be small, and it is at least as large as \(\epsilon_4(0)^2\). In our notation, \(\epsilon_4(0)^2=\dimp^2\|\nabla^4\lgd(\map)\|_{\DFLG}\). Using the arguments in Appendix~\ref{app:sub:third}, it is straightforward to show that if \(h^{(4)}\) is of constant sign in the interval \([-\|\calR q_{\map}\|_\infty, \|\calR q_{\map}\|_\infty]\), then \(\epsilon_4(0)^2\gtrsim (m\wedge\dimp)\dimp^2/n\). Therefore, the bound of~\cite{katskew} requires at least \(n\gg\dimp^2\), so it is not dimension-free. For this reason, and due to the different structure of the bound, we do not further investigate this bound in the present work. 
\end{remark}

\begin{table}
\renewcommand{\arraystretch}{1.3}
\begin{tabular}{c|c|c|c|c|c}
 &  $\DFLG$ &  $I_{\dimp}/\normG(I_{\dimp})$ & $\DFL(\gamma_0^*)$ & $\mathrm{UB}(\DFLG)/\mathrm{UB}(\DFL(\gamma_0^*))$ & $\mathrm{UB}(I_{\dimp}/\normG(I_{\dimp}))/\mathrm{UB}(\DFL(\gamma_0^*))$\\
\hline
\(\dimp<m\) & \(\dimp^3/n\) &  \(\dimp^{2+6\beta}/n\)&  \(\dimp^3/n\)   & 1 & \(\dimp^{6\beta-1}\)\\
\hline
\(m<\dimp<m_0^*\)& \(\dimp^2m/n\) &  \(m^{2+6\beta}/n\)& \( m^2\dimp/n\) & \(\dimp/m\) & \(m^{6\beta}/\dimp\)\\
\hline
\(m_0^*<\dimp\) & \(\dimp^2m/n\)&  \(m^{2+6\beta}/n\) & \(m^2m_0^*/n\)  & \((\dimp/m)(\dimp/m_0^*)\) & \(m^{6\beta}/m_0^*\)\\
\end{tabular}
\caption{The second, third, and fourth columns display the squared local TV upper bound, \(\mathrm{UB}(\DFL):=\dimL(\DFL)^2\mathrm{UB}(\|\nabla^3\lgd\|_{\DFL,\infty})^2\) for three choices of \(\DFL\). $\DFLG$ is used in~\cite{katsBVM}, $I_{\dimp}/\normG(I_{\dimp})$ is used in~\cite{bp}, and $\DFL(\gamma_0^*)$ is our optimized choice. The last two columns show the factor by which the $\DFL(\gamma_0^*)$-based bound improves on the other two bounds, respectively. In the setting considered in Lemma~\ref{lma:compare}, where \(2\gamma+2\beta>2\) and \(\beta>1/2\), each of the ratios in the last two columns is greater than 1. Thus \(\DFL(\gamma_0^*)\) improves on both \(\DFLG\) and \(I_{\dimp}/\normG(I_{\dimp})\).}
\label{table:cases}
\end{table}

\subsection{Additional details on the choice of \(\DFL\)}\label{subsec:gam0}
In this section, we explain the rationale to consider the family \(\DFL^2(\gamma_0)\) given in~\eqref{DFLgam0def}. We then show that \(\DFL(\gamma_0)\) is close to a diagonal matrix, and use this to compute \(\dimL(\DFL(\gamma_0))\) and \(\|\nabla^3\lgd\|_{\DFL(\gamma_0)}\). Finally, we optimize over \(\gamma_0\) to arrive at \(\gamma_0^*=\gamma-1/2-1/2m\). 

We are interested in applying Corollary~\ref{corr:DD0} to bound the TV error of the LA to \(\PfL\), with \(\lgd\) as in~\eqref{lgd:inv}. 
Our setup is exactly as in Section~\ref{subsub:struct}. Thus, as discussed in Section~\ref{Sconc}, to ensure \(\nabla^2\lgd\succeq\DFLG^2-\DFL^2\), we should take \(\DFL^2\) of the form \(\DFL^2=\nabla^2L(\map)+ G_0^2\) where \(0\preceq G_0^2\preceq G^2\). Since \(G^2\) is diagonal with the $k$th diagonal entry growing as a power of \(k\), it is natural to take \(G_0^2\) of the same form. This leads to
$$\DFL^2(\gamma_0)=\nabla^2L(\map)+G_0^2,\qquad G_0^2=\diag(k^{2\gamma_0})_{k=1}^{\dimp}, \qquad\gamma_0\leq\gamma.$$ The condition \(\normG(\DFL(\gamma_0))=1\) is satisfied because \(\DFLG^2\succeq\DFL^2(\gamma_0)\) and \(\|\DFL(\gamma_0)\ev_1\|=\|\DFLG\ev_1\|\), where \(\ev_1=(1,0,\dots,0)\). Note also that \(\DFL^2(-\infty)=\nabla^2L(\map)+\diag(1,0,\dots,0)\) and \(\DFL^2(\gamma)=\nabla^2L(\map)+G^2=\DFLG^2\). Thus by letting \(\gamma_0\) range between \(-\infty\) and \(\gamma\), the matrix \(\DFL^2(\gamma_0)\) ranges from (almost) the Hessian of the negative log likelihood, to the Hessian of the negative log posterior. 

\begin{remark}As we can see, the primary rationale for choosing this family is simply by necessity: to take advantage of the concentration result Theorem~\ref{thm:conc}, we need to satisfy~\eqref{n2lgd}, and this restricts us to matrices \(\DFL^2=\nabla^2L(\map)+G_0^2\). The concentration result is what gives the exponentially small second term in the TV bound of Corollary~\ref{corr:DD0}. Thus it is interesting that, to find \(\DFL^2\) that leads to a small \emph{local} contribution to our TV bound, it turn out to be a good idea to impose the assumption needed for a very small \emph{tail} contribution.\end{remark}

Next, we state diagonal upper and lower bounds on \(\DFL^2(\gamma_0)\). 
\begin{lemma}\label{lma:hess} 
Under Assumption~\ref{assume:Rqinfty}, it holds
\begin{EQA}\label{ch-Ch-prelim}
&&\DFL^2(\gamma_0)\asymp R^\T R + G_0^2.
\end{EQA}
If Assumption~\ref{assume:ortho} is also satisfied, then for all \(\gamma_0\) such that \(2\gamma_0+2\beta\geq\lambda\), 
\begin{EQA}\label{eq:hess}
&&\DFL^2(\gamma_0)\asymp \diag\left(\frac{n}{k^{2\beta}}+k^{2\gamma_0}\right)_{k=1}^{\dimp}.
\end{EQA} 
\end{lemma}
To see why Assumption~\ref{assume:Rqinfty} is needed, recall from~\eqref{LL-hess-radon} the form of \(\nabla^2L(\map)\), and recall that \(R_j^\T\map=(\calR q_{\map})(j/n)\). If \((\calR q_{\map})(x)\) goes off to \(\pm\infty\), then the factors \(h''(R_j^\T\map)\) in the summands making up \(\nabla^2L(\map)\) may be arbitrarily close to zero. This is because \(h\) is convex but typically not strongly convex.

The diagonal bounds~\eqref{eq:hess} significantly simplify the calculation of \(\dimL(\DFL)\) and \(\|\nabla^3f\|_{\DFL}\). We now show that these two quantities are equivalent to, and bounded by, the sums \(S_{\dimL}\) and \(S_{\dltwu}\), respectively. Recall the definition of these sums from~\eqref{SdimL-def} and~\eqref{Stau-def}. 
\begin{lemma}\label{lma:effdim:radon}Suppose Assumptions~\ref{assume:bdd},~\ref{assume:ortho},~\ref{assume:Rqinfty} are all satisfied, and \(\gamma_0\) is as in Lemma~\ref{lma:hess}. Then \(\dimL(\DFL(\gamma_0))\asymp S_{\dimL}\), and
\begin{EQA}
&&\sup_{\|\DFL(\gamma_0)\uv\|\leq\rrdc}\|\nabla^3\lgd(\map+\uv)\|_{\DFL(\gamma_0)}\les  D_3(C(1\vee \rrdc S_{\dltwu}))S_{\dltwu},\label{tau3bound-0}
\end{EQA} 
where \(D_3(K):=\sup\left\{|h'''(t)|\; : \; |t|\leq  K\right\}\). 
\end{lemma}
The bound~\eqref{tau3bound-0} is the result for which we need Assumption~\ref{assume:bdd}. Specifically, proving~\eqref{tau3bound-0} involves bounding \(\sup_{\|\DFL(\gamma_0)\uv\|=1}\|\calR q_{\uv}\|_\infty\), and we use that \(\|\calR q_{\uv}\|_\infty\leq \sum_{k=1}^{\dimp}|u_k|\sqrt{\lambda_k}\|\psi_k\|_\infty\). From here, we use that \(\sup_k\|\psi_k\|_\infty\les 1\). See Appendix~\ref{app:sub:third} for more details.

Lemma~\ref{lma:effdim:radon} and Corollary~\ref{corr:DD0} with \(\rrdc=1/\sqrt{S_{\dltwu}}\) now give Theorem~\ref{corr:TV:radon} (see Appendix~\ref{app:sec3:proofs} for the full proof). We recall that the theorem shows \(\TV(\PfL,\lap)\les S_{\dimp}S_{\dltwu}\) if \(S_{\dimp}S_{\dltwu}\) is small enough and \(2\gamma_0+2\beta\geq\lambda\). We now optimize the upper bound \(S_{\dimp}S_{\dltwu}\) over \(\gamma_0\). To see how to choose \(\gamma_0\), it is instructive to consider the case \(\dimp\to\infty\). Let \begin{EQA}\label{def:m-m0}
&&m=n^{1/(2\beta+2\gamma)},\qquad m_0=n^{1/(2\beta+2\gamma_0)}\geq m.
\end{EQA} For \(k<m\), the denominator in the summands of \(S_{\dimL}\) is dominated by \(n\), and for \(k>m\) the denominator is dominated by \(k^{2\beta+2\gamma}\). The value \(m_0\) is the analogous threshold for the denominator of the summands of \(S_{\dltwu}\). It is then clear that as \(\dimp\to\infty\), we have
\begin{EQA}
S_{\dimL}(\gamma_0)&\asymp &m + \sum_{k>m}nk^{-2\gamma-2\beta}+\sum_{k>m}k^{2\gamma_0-2\gamma}\asymp m(1+m^{2\gamma_0-2\gamma})\asymp m,\label{Sdimp}\\
S_{\dltwu}(\gamma_0)^2&\asymp &\frac{m_0}{n} + \sum_{k>m_0}k^{-2\beta-2\gamma_0}\asymp \frac{m_0}{n}+m_0^{1-2\beta-2\gamma_0}\asymp\frac{m_0}{n},
\end{EQA} where the first bound holds if \(\gamma_0\) is bounded above away from \(\gamma-1/2\) and the second bound holds if \(\gamma_0\) is bounded below away from \(1/2-\beta\). But now, the optimal choice of \(\gamma_0\) is clear. The value of \(S_{\dimL}\) does not depend on \(\gamma_0\) at all when \(\gamma_0\) is bounded above away from \(\gamma-1/2\). Meanwhile, increasing \(\gamma_0\) decreases \(m_0\) and hence decreases \(S_{\dltwu}\). Thus we should take \(\gamma_0\) as close as possible to \(\gamma-1/2\). In Appendix~\ref{app:Sbounds}, we determine how close \(\gamma_0\) can get to \(\gamma-1/2\) without the constant of proportionality blowing up. We obtain the optimal choice is \(\gamma_0^*=\gamma-1/2-1/2m\), and the following bounds.
\begin{lemma}\label{lma:Sranges}Suppose \(2\beta+2\gamma>2\) and \(n>(\beta+\gamma-1)^{-2\beta-2\gamma}\). Then
\begin{EQA}
&&S_{\dimL}(\gamma_0^*)\asymp m\wedge\dimp,\qquad S_{\dltwu}(\gamma_0^*)\asymp\frac{m_0^*\wedge\dimp}{n},
\end{EQA}
where \(\gamma_0^*=\gamma-1/2-1/2m\) and \(m_0^*=n^{1/(2\beta+2\gamma_0^*)}\).
\end{lemma}
\begin{remark}\label{rk:effd}A salient feature in the optimization is that the effective dimension \(\dimL(\DFL)\) is more stable to changes in \(\DFL\) than is \(\dltwu_3(\UV,\DFL)\). For example, within the family \(\DFL(\gamma_0)\), the value \(\dimL(\DFL(\gamma_0))\) stays the same for all \(\gamma_0\) bounded above away from \(\gamma-1/2\); see~\eqref{Sdimp}, and recall that \(\dimL(\DFL(\gamma_0))\asymp S_{\dimL}(\gamma_0)\). The very different choice \(\DFL=I_{\dimp}/\normG(I_{\dimp})\) also leads to the same \(\dimL(\DFL)\), as seen in Lemma~\ref{lma:compare}. This single value of \(\dimL(\DFL)\) arising for the various choices of \(\DFL\) is \(\dimL(\DFL)=m=n^{1/(2\beta+2\gamma)}\) (for large values of \(\dimp\)), which has the natural interpretation of being the number of coordinates for which the information in the likelihood outweighs the prior regularization. (Recall that \(\nabla^2L(\map)\approx \diag(n/k^{2\beta})\).) 
\end{remark}

\appendix
%
\section{Proof of concentration result}\label{app:conc}


In this section, we use the shorthands \(\UV=\UV(\DFc,\rrdc)\), \(\dltwb=\dltwb(\UV,\DFc)\), \(\dimL=\dimL(\DFc)\). We also redefine \(\UV=\UV-\map\) in this section only. The proof of Theorem~\ref{thm:conc} follows from the following two lemmas, which are proved at the end of the section.
\begin{lemma}\label{lma:concentrate:1}
Under the conditions of Theorem~\ref{thm:conc}, it holds
\begin{EQA}
	S\eqdef \frac{\int_{\UV^c} \ex^{f(\map)-f(\map+\uv)} \, d\uv}
		 {\int \ex^{ - \| \DFLG \uv \|^{2}/2} \, d \uv}	& \leq &
	\frac13\ex^{-\dimL/2}\ex^{-\frac13(\rrdc-3\sqrt{\dimL})^2}.
\label{fiinIHUniUef}
\end{EQA} 
\end{lemma}
\begin{lemma}\label{lma:concentrate:2}
Under the conditions of Theorem~\ref{thm:conc}, it holds
\begin{EQA}\label{bulkint}
	\int_{\UV} \ex^{\lgd(\map)-\lgd(\map+\uv)} \, d\uv  
	& \geq &
	\ex^{ - \dltwb \, \dimL/2} \, \int_{\UV} \ex^{- \| \DFLG \uv \|^{2}/2} d\uv \, .
\label{emdw3dp0iUfxu}
\end{EQA} 
\end{lemma} 
\begin{proof}[Proof of Theorem~\ref{thm:conc}]
Let \( s=\ex^{-(\rrdc -\sqrt{\dimL})^2/2} \)  and \( \lgd(\map;\uv)=\lgd(\map+\uv)-\lgd(\map) \).  Then
\begin{EQA}
 \int_{\UV} \ex^{-\lgd(\map;\uv)} \, d\uv  & \geq & \ex^{ - \dltwb \, \dimL/2} \, \int_{\UV} \ex^{- \| \DFLG \uv \|^{2}/2} d\uv\\
& \geq & \ex^{ - \dltwb \, \dimL/2}(1-s)\int\ex^{- \| \DFLG \uv \|^{2}/2} d\uv 
\end{EQA} by Lemma~\ref{lma:concentrate:2} and~\eqref{gauss-tail-0}. Using this inequality, we have
\begin{EQA}
\PfL((\map+\UV)^c) &=& \frac{\int_{\UV^{c}} \ex^{-\lgd(\map;\uv)} \, d\uv}{\int  \ex^{-\lgd(\map;\uv)} \, d\uv}= \frac{\int_{\UV^{c}} \ex^{-\lgd(\map;\uv)} \, d\uv}{\int_{\UV} \ex^{-\lgd(\map;\uv)} \, d\uv + \int_{\UV^{c}}\ex^{-\lgd(\map;\uv)} \, d\uv}\\
&\leq & \frac{\int_{\UV^{c}} \ex^{-\lgd(\map;\uv)} \, d\uv}{ \ex^{ - \dltwb \, \dimL/2}(1-s)\int\ex^{- \| \DFLG \uv \|^{2}/2} d\uv + \int_{\UV^{c}}\ex^{-\lgd(\map;\uv)} \, d\uv}
\end{EQA}
Now, note that $a\mapsto a/(b+a)$ is increasing. Thus we can use that
\begin{EQA}
\int_{\UV^{c}} \ex^{f(\map;\uv)} \, d\uv &\leq &t \int\ex^{ - \| \DFLG \uv \|^{2}/2} \, d \uv, \qquad t:=\frac13\ex^{-\dimL/2}\ex^{-\frac13(\rrdc-3\sqrt{\dimL})^2}
\end{EQA} by Lemma~\ref{lma:concentrate:1}. Combining the above two inequalities we conclude that
\begin{EQA}
\PfL((\map+\UV)^c) &\leq & \frac{t }{\ex^{-\dltwb \,\dimL/2}(1-s) +t} \leq \frac{t }{\ex^{-\dltwb \,\dimL/2} +(t-s)}\leq \ex^{\dltwb \,\dimL/2}t.
\end{EQA} The last inequality uses that $t\geq s$, which is straightforward to check.  Since \(\dltwb\leq1/3\) we have \(\ex^{\dltwb \,\dimL/2}t \leq\frac13\ex^{-\dimL/3- (\rrdc - 3\sqrt{\dimL})^2/3}\), concluding the proof. 
\end{proof}
To prove Lemma~\ref{lma:concentrate:1}, we give two auxiliary lemmas.
\begin{lemma}\label{lma:lambda}
Let  \(\lgd_0(\uv):=\lgd(\map+\uv)-\lgd(\map)-\frac12\uv^{\T}(\DFLG^2-\DFc^2)\uv\). It holds
\begin{EQA}
\lgd_0(\uv)&\geq& \frac{1-\dltwb}{2}\rrdc\|\DFc\uv\|,\qquad\forall\,\|\DFc\uv\|\leq\rrdc.\label{lgd-dltwb}\\
\end{EQA}
\end{lemma}
\begin{remark}\label{rk:conc1}The proof of Proposition A.12 in~\cite{spok23} contains a result analogous to this lemma. However, the second inequality in (A.24) does not follow from the definition of \(\dltwb\).We correct this in the below proof.
\end{remark}
\begin{proof}
Note that \(\lgd_0\) is convex by~\eqref{n2lgd}. Fix \(\uv\) such that \(\|\DFc\uv\|\geq\rrdc\), and let \(\uv_0\) be a vector in the same direction as \(\uv\) but with \(\|\DFc\uv_0\|=\rrdc\). By convexity, we have 
\begin{EQA}\label{cvx1}
&&\lgd_0(\uv)-\lgd_0(0)\geq \frac{\|\DFc\uv\|}{\rrdc}\left(\lgd_0(\uv_0)-\lgd_0(0)\right).
\end{EQA} 
Next, since \(\lgd_0\) differs from \(\lgd\) by a quadratic function, we have
\begin{EQA}
\bigl|\lgd_0(\uv_0) &-& \lgd_0(0)-\frac12\langle\nabla^2\lgd_0(0), \uv_0^{\otimes2}\rangle\bigr| \\
&=&\bigl|\lgd(\map+\uv_0) - \lgd(\map)-\frac12\langle\nabla^2\lgd(\map), \uv_0^{\otimes2}\rangle\bigr| \\
&=& |\dltw_3(\uv_0)| \leq \frac{\dltwb}{2}\|\DFc\uv_0\|^2 = \frac{\dltwb}{2}{\rrdc}^2,
\end{EQA} by the definition~\eqref{def:omega} of \(\dltwb\). Using that \(\langle\nabla^2\lgd_0(0), \uv_0^{\otimes2}\rangle = \|\DFc\uv_0\|^2={\rrdc}^2\), the above display gives
\begin{EQA}
&&\lgd_0(\uv_0) - \lgd_0(0) \geq \frac{1-\dltwb}{2}{\rrdc}^2.
\end{EQA} Combining this with~\eqref{cvx1} concludes the proof.
\end{proof}
\begin{lemma}\label{lma:R}Define
\begin{EQA}
R(\lambda) &=& \E\left[\exp \bigl(- \lambda \| \Tau \gaussv \| + \| \Tau \gaussv \|^{2} / 2 \bigr)\Ind\bigl( \| \Tau \gaussv \| > \rrdc \bigr) \right],
\end{EQA} where \(\Tau=\DFc\DFLG^{-1}\) and \(Z\sim\ND(0, I_{\dimp})\). If \(\rrdc>\sqrt{\dimL}\), then
\begin{EQA}
R(\lambda) &\leq & \ex^{-\dimL/2}\ex^{-b\rrdc }\left(1+\frac{1}{b^2}+\frac{\rrdc-\lambda}{b}\right),
\end{EQA} where \(b=\lambda-\sqrt{\dimL}\).
\end{lemma}
\begin{proof}
We can write \(R\) as the following one-dimensional integral: 
\begin{EQA}\label{Rrr-0}
	R&=&\int_{\rrdc}^{\infty}\ex^{-\lambda\zq +\frac{\zq^2}{2}} \, 
		d\P\bigl( \| \Tau \gaussv \| > \zq \bigr)	\\
	&=& -\ex^{-\lambda\zq +\frac{\zq^2}{2}}\P\bigl( \| \Tau \gaussv \| >\zq \bigr)\bigg\vert_{\zq=\rrdc}^{\zq=\infty} + \int_{\rrdc}^{\infty}(\zq-\lambda)\ex^{-\lambda\zq+\frac{\zq^2}{2}}\P\bigl( \| \Tau \gaussv \| > \zq \bigr)d\zq\,.\\
	&\leq &
	\ex^{-\lambda\rrdc +\frac{{\rrdc}^2}{2}}\P\bigl( \| \Tau \gaussv \| >\rrdc \bigr) + \int_{\rrdc}^{\infty}(\zq-\lambda)\ex^{-\lambda\zq+\frac{\zq^2}{2}}\P\bigl( \| \Tau \gaussv \| > \zq \bigr)d\zq\,.
\end{EQA}
Now, using Corollary~\ref{corr:gaussconc} on Gaussian concentration, we have for any \( \zq \geq \sqrt{\dimL} \),
\begin{EQA}
	\P\bigl( \| \Tau \gaussv \| > \zq \bigr)
	& \leq & 
	\exp\bigl( - \frac12(\zq - \sqrt{\dimL})^{2} \bigr).
\label{2emxPTgasps2d}
\end{EQA}

Therefore,
\begin{EQA}\label{exlam}
\ex^{-\lambda\zq +\frac{\zq^2}{2}}\P\bigl( \| \Tau \gaussv \| >\zq \bigr)&\leq&\exp\bigl(-\lambda\zq +\frac{\zq^2}{2} - (\zq - \sqrt{\dimL})^{2}/2 \bigr)\\
&=&\ex^{-\dimL/2}\exp\bigl(-(\lambda-\sqrt{\dimL})\zq\bigr)
\end{EQA} Substituting this bound into the last line of~\eqref{Rrr-0} gives
\begin{EQA}\label{Rrr}
R&\leq& \ex^{-\dimL/2}\ex^{-(\lambda-\sqrt{\dimL})\rrdc }+ \ex^{-\dimL/2}\int_{\rrdc}^{\infty}(\zq-\lambda)\ex^{-(\lambda-\sqrt{\dimL})\zq}d\zq\\
&=&\ex^{-\dimL/2}\ex^{-b\rrdc }+ \ex^{-\dimL/2}\int_{\rrdc}^{\infty}(\zq-\lambda)\ex^{-b\zq}d\zq\\
&=& \ex^{-\dimL/2}\ex^{-b\rrdc }\left(1+\frac{1}{b^2}+\frac{\rrdc-\lambda}{b}\right).
\end{EQA}
\end{proof}
\begin{proof}[Proof of Lemma~\ref{lma:concentrate:1}]
Let \(X\sim\mathcal N(0, \DFLG^{-2})\) and \(\gaussv\sim\mathcal N(0, I_{\dimp})\). We can write \(S\) as
\begin{EQA}\label{S1}
S &=&\frac{\int  \Ind\bigl( \| \DFc \uv \| > \rrdc \bigr)\exp \bigl( -\lgd_0(\uv)+ \| \DFc \uv \|^{2}/{2} - \| \DFLG \uv \|^{2}/{2}\bigr) \, d \uv}{\int \exp \bigl( - \| \DFLG \uv \|^{2}/{2} \bigr) \, d \uv}\\
&=&\E\left[\exp \bigl(- \lgd_0(X) +\|\DFc X\|^2/2\bigr)\Ind\bigl( \| \DFc X \| > \rrdc \bigr) \right].
\end{EQA}
Next, let \(\lambda = (1-\dltwb)\rrdc/2\). Using Lemma~\ref{lma:lambda}, we have
\begin{EQA}
S&\leq& \E\left[\exp \bigl(- \lambda\|\DFc X\| +\|\DFc X\|^2/2\bigr)\Ind\bigl( \| \DFc X \| > \rrdc \bigr) \right]\\
&=&\E\left[\exp \bigl(	- \lambda \| \Tau \gaussv \| + \| \Tau \gaussv \|^{2} / 2 \bigr)\Ind\bigl( \| \Tau \gaussv \| > \rrdc \bigr) \right]=R(\lambda),
\end{EQA} where \(\Tau=\DFc\DFLG^{-1}\), and \(R\) is as in Lemma~\ref{lma:R}. From here, we use Lemma~\ref{lma:R} and the definition of \(\lambda\) to get
\begin{EQA}
S&\leq & R(\lambda) \leq \ex^{-\dimL/2}\ex^{-(\lambda-\sqrt{\dimL})\rrdc }\left(1+\frac{1}{(\lambda-\sqrt{\dimL})^2}+\frac{1+\dltwb}{2}\frac{\rrdc}{\lambda-\sqrt{\dimL}}\right)
\end{EQA}Using the assumption \(\dltwb\leq1/3\) gives \(\lambda\geq\rrdc/3\), and hence the above further simplifies as
\begin{EQA}
S&\leq & R(\lambda) \leq \ex^{-\dimL/2}\ex^{-(\rrdc/3-\sqrt{\dimL})\rrdc }\left(1+\frac{1}{(\rrdc/3-\sqrt{\dimL})^2}+\frac{2}{3}\,\frac{\rrdc}{\rrdc/3-\sqrt{\dimL}}\right)
\end{EQA}
Next, we use that \(\rrdc/3-\sqrt{\dimL}\geq1\) by assumption to get
\begin{EQA}
1+\frac{1}{(\rrdc/3-\sqrt{\dimL})^2}+\frac{2}{3}\,\frac{\rrdc}{\rrdc/3-\sqrt{\dimL}} &=& 3+\frac{1}{(\rrdc/3-\sqrt{\dimL})^2} + \frac{2\sqrt{\dimL}}{\rrdc/3-\sqrt{\dimL}}\\
&\leq& 4+2\sqrt{\dimL}.
\end{EQA}
Furthermore,
\begin{EQA}
\frac{\dimL}{2} + (\rrdc/3-\sqrt{\dimL})\rrdc &=& \frac13(\rrdc-3\sqrt{\dimL})^2+ \rrdc\sqrt{\dimL}-2.5\dimL\\
&=&\frac13(\rrdc-3\sqrt{\dimL})^2+ (\rrdc-3\sqrt{\dimL}-3)\sqrt{\dimL} + \frac{\dimL}{2} + 3\sqrt{\dimL}\\
&\geq&\frac13(\rrdc-3\sqrt{\dimL})^2 +  \frac{\dimL}{2} + 3\sqrt{\dimL}.
\end{EQA}
Combining the last three displays gives \(S\leq (4+2\sqrt{\dimL})\ex^{-\dimL/2-3\sqrt{\dimL}}\ex^{-\frac13(\rrdc-3\sqrt{\dimL})^2}\). Finally, we use that \((4+2\sqrt{\dimL})\exp(-3\sqrt{\dimL})\leq 1/3\) for all \(\dimL\geq1\) to conclude.
\end{proof}

\begin{remark}\label{rk:conc2}Proposition A.6 in~\cite{spok23} is analogous to our Lemma~\ref{lma:concentrate:2}. However, the change of variables described below (A.14) does not yield the upper bound in the second line of the following display. We correct this in the below proof. \end{remark}

\begin{proof}[Proof of Lemma~\ref{lma:concentrate:2}]
By Definition~\ref{def:del3-omega} of \( \dltw_3 \) and \(\dltwb\), we have
\begin{EQA}
	\int_{\UV} \ex^{\lgd(\map)-\lgd(\map+\uv)} \, d\uv
	& = &
	\int_{\UV} \ex^{- \| \DFLG \uv \|^{2}/2 + \dltw_{3}(\uv)} \, d\uv \\
	&\geq &
	\int_{\UV} \ex^{- \| \DFLG \uv \|^{2}/2 - \dltwb \| \DFc\uv \|^{2}/2} \, d\uv =\int_{\UV} \ex^{- \frac12\uv^{\T}A\uv} d\uv,
\label{PhtiUexmHu222}
\end{EQA}where $A=\DFLG^2+\dltwb\DFc^2$.
Now, let $B^2= \Id + \dltwb \, \DFLG^{-1} \DFc^{2} \, \DFLG^{-1}$, and consider the change of variables \(\uv=\DFLG^{-1} B^{-1}\DFLG\wv$. Then
$$
\uv^{\T}A\uv = \wv^{\T}\DFLG B^{-1}\DFLG^{-1}A\DFLG^{-1} B^{-1}\DFLG\wv= \wv^{\T}\DFLG^2 \wv,
$$using that $\DFLG^{-1}A\DFLG^{-1}=B^2$. Also note that $\det(\DFLG^{-1} B^{-1}\DFLG)=\det(B^{-1})$. Hence the change of variables gives
\begin{EQA}
	\int_{\UV} \ex^{\lgd(\map)-\lgd(\map+\uv)} \, d\uv
	& \geq &
	\det \bigl( \Id + \dltwb \, \DFLG^{-1} \DFc^{2} \, \DFLG^{-1} \bigr)^{-1/2} 
	\int_{\DFLG^{-1} B\DFLG\UV} \ex^{- \| \DFLG \wv \|^{2}/2 } \, d\wv \\
	&\geq & \ex^{-\frac{\dltwb\dimL}{2}}\int_{\DFLG^{-1} B\DFLG\UV} \ex^{- \| \DFLG \wv \|^{2}/2 } \, d\wv.
\label{iUefxudetId3H121}
\end{EQA}
To get the second inequality we used that $\log(1+x)\leq x$, $x\geq0$ to deduce that
\begin{EQA}
\label{12Ip2t3t4Hp3}
	\log \det\bigl( \Id + \dltwb\DFLG^{-1} \DFc^{2} \, \DFLG^{-1} \bigr)
	& \leq &
	\dltwb \, \tr(\DFLG^{-1} \DFc^{2} \, \DFLG^{-1}) =\dltwb\dimL.
\label{12Im2t3t4Hm30vfufdkjr}
\end{EQA}
To conclude the proof, we show $\DFLG^{-1} B\DFLG\UV\supset\UV$ which is equivalent to showing $\UV\supset \DFLG^{-1} B^{-1}\DFLG\UV$. Now, an element of the set $\DFLG^{-1} B^{-1}\DFLG\UV$ is given by $\vv=\DFLG^{-1} B^{-1}\DFLG\uv$ for $\uv\in\UV$, so that $\|\DFc\uv\|\leq\rrdc$. We need to show that $\vv\in\UV$, i.e. that $\|\DFc\vv\|\leq \rrdc$. But 
\begin{EQA}
\|\DFc\vv\|& =& \|\DFc\DFLG^{-1} B^{-1}\DFLG\uv\| \leq \|\DFc\DFLG^{-1} B^{-1}\DFLG\DFc^{-1}\|\|\DFc\uv\| \leq \|\DFc\DFLG^{-1} B^{-1}\DFLG\DFc^{-1}\|\rrdc
\end{EQA}
Thus it remains to show $ \|S B^{-1}S^{-1}\|\leq1$, where $S= \DFc\DFLG^{-1}$. Let $S=U\Lambda V^{\T}$ be the singular value decomposition of $S$. Then
\begin{EQA}
B^2&=&\Id+\dltwb S^{\T}S=V(\Id+\dltwb \Lambda^2)V^{\T}
\end{EQA} and hence $B =V(\Id+\dltwb \Lambda^2)^{1/2}V^{\T}$. Therefore,
\begin{EQA}
SB^{-1}S^{-1} &=& U\Lambda V^{\T}V(\Id+\dltwb\Lambda^2)^{-1/2}V^{\T}V\Lambda^{-1} U^{\T} = U(\Id+\dltwb \Lambda^2)^{-1/2}U^{\T},
\end{EQA} from which it is now clear that \(\| SB^{-1}S^{-1}\|\leq 1 \). 
\end{proof}

\section{Details of comparison to prior work}\label{app:compare}

As noted in Section~\ref{sec:comparegen} of the main text, \cite{katsBVM} also starts from the decomposition~\eqref{TVdecomp}, and then further bounds both the local TV distance and the two tail probabilities. 
The result is Theorem 2.2, which takes the form
\begin{EQA}\label{TVkatsBVM}
&&\TV(\PfL,\lap)\leq \frac{\delta_3(r)\dimp}{\sqrt{n}} +  \mathrm{UB}(\PfL(\UV^c)) + \mathrm{UB}(\lap(\UV^c)),
\end{EQA} where \(\UV=\{\|\DFLG\uv\|\leq r\sqrt{\dimp}\}\). The local bound requires \(r\delta_3(r)\sqrt{\dimp/n}\leq1/2\) and the tail bound requires convexity of \(\lgd\). The quantity \(\delta_3\) is defined in (2.1). It is straightforward to show that \(\delta_3(r)/\sqrt n =\dltwu_3(\UV,\DFLG)\), and furthermore, we have \(\dimp=\Tr(\DFLG^{-2}\DFLG^2)=\dimL(\DFLG)\). Finally, \(r\delta_3(r)\sqrt{\dimp/n}\leq1/2\) implies \(\dltwb_3(\UV,\DFLG)\leq 1/2\), as in our Corollary~\ref{thm:TV}. Thus under a slightly stronger condition than our Corollary~\ref{thm:TV},~\cite{katsBVM} obtains the result stated in~\eqref{katsBVM-simple-ub}.

Next, we compare our results to Theorem 17 in~\cite{bp}. The following is the key notation in this work. 
\begin{equation}\label{translate-bp}\begin{gathered}
\UV=\{\|\thv-\map\|\leq\bar\delta\},\\
 \bar J_n(\bar\theta_n) = n^{-1}\DFLG^2,\quad \Tr(\bar J_n(\bar\theta_n)^{-1})=n\Tr(\DFLG^{-2}),\\
\bar M_2 = n^{-1}\sup_{\thv\in\UV}\|\nabla^3f(\thv)\|,\\
\bar\lambda_{\min}(\bar\theta_n)=n^{-1}\lambda_{\min}(\DFLG^2) = n^{-1}\|\DFLG^{-1}\|^{-2} = n^{-1}\normG(I_{\dimp})^{-2},\\
\bar{\mathcal D}(n,\bar\delta) \geq \lap(\UV^c).
\end{gathered}\end{equation}
The fact that \(\bar{\mathcal D}(n,\bar\delta)\) is an upper bound on \(\lap(\UV^c)\) follows by the definition of \(\bar{\mathcal D}(n,\bar\delta)\) on page 13, and the standard Gaussian tail bound~\eqref{gauss-tail-0}.
Throwing out constant factors, the authors' bound in Theorem 17 is the following:
\begin{EQA}\label{bp-simple-ub-app}
&&\TV(\PfL,\lap)\les \frac{\Tr(\bar J_n(\bar\theta_n)^{-1})\bar M_2}{\sqrt{n(\bar\lambda_{\min}(\bar\theta_n)-\bar\delta\bar M_2)(1-\bar{\mathcal D}(n,\bar\delta))}} +  A_2n^{d/2}e^{-n\bar\kappa} + \bar{\mathcal D}(n,\bar\delta)
\end{EQA}
Now, the last summand, \( \bar{\mathcal D}(n,\bar\delta)\), is an upper bound on \(\lap(\UV^c)\). The second summand, \(A_2n^{d/2}e^{-n\bar\kappa}\), can be shown to be an upper bound on \(\PfL(\UV^c)\) under the authors' assumptions. This follows by the six lines of calculations on page 48 of~\cite{bp}, below the definition of \(I_{2,1}^{MAP}\) and \(I_{2,2}^{MAP}\). Thus it remains to show the first summand on the righthand side of~\eqref{bp-simple-ub-app} can be brought into the form of the first summand on the righthand side of~\eqref{bp-simple-ub}. We immediately have that  \(1- \bar{\mathcal D}(n,\bar\delta)\) is a lower bound on \(\lap(\UV)\). Next, note that
\begin{EQA}
n(\bar\lambda_{\min}(\bar\theta_n)-\bar\delta\bar M_2)&=&\normG(I_{\dimp})^{-2}- \sup_{\thv\in\UV}\|\nabla^3\lgd(\thv)\|_{I_{\dimp}}\sup_{\thv\in\UV}\|\thv-\map\|\\
&\leq&\normG(I_{\dimp})^{-2}-\dltwb_3(\UV,I_{\dimp})= \lambda(\UV,I_{\dimp}),
\end{EQA} using~\eqref{del3bds2}. Thus the denominator of the first term in~\eqref{bp-simple-ub-app} can indeed be expressed as $\sqrt{\mathrm{LB}(\lambda(\UV,I_{\dimp}))\mathrm{LB}(\lap(\UV))}$. Finally, using the second and third lines in~\eqref{translate-bp}, we see that the numerator of the first term in~\eqref{bp-simple-ub-app} is precisely \(\sup_{\thv\in\UV}\|\nabla^3\lgd(\thv)\|_{I_{\dimp}}\Tr(\DFLG^{-2})\). We have therefore brought~\eqref{bp-simple-ub-app} into the form~\eqref{bp-simple-ub}. Next, we show how we can normalize this bound to bring it into the form~\eqref{bp-ub-norm}. Define 
$$
\mathrm{LB}(\lambda(\UV,\DFL))=\normG(\DFL)^{-2}- \sup_{\thv\in\UV}\|\nabla^3\lgd(\thv)\|_{\DFL}\sup_{\thv\in\UV}\|\DFL(\thv-\map)\|.$$
Then it is clear that
\begin{EQA}
\|\nabla^3\lgd(\thv)\|_{I_{\dimp}}&=&\|\nabla^3\lgd(\thv)\|_{I_{\dimp}/\normG(I_{\dimp})}\,\normG(I_{\dimp})^{-3},\\
\Tr(\DFLG^{-2})&=&\dimL(I_{\dimp}/\normG(I_{\dimp}))\,\normG(I_{\dimp})^2,\\
\mathrm{LB}(\lambda(\UV,I_{\dimp}))&=&\mathrm{LB}(\lambda(\UV,I_{\dimp}/\normG(I_{\dimp})))\normG(I_{\dimp})^{-2}.
\end{EQA}
Substituting these expressions into the first term of~\eqref{bp-simple-ub} gives
$$
\frac{\sup_{\thv\in\UV}\|\nabla^3\lgd(\thv)\|_{I_{\dimp}/\normG(I_{\dimp})}\,\dimL(I_{\dimp}/\normG(I_{\dimp}))}{\sqrt{\mathrm{LB}(\lambda(\UV,I_{\dimp}/\normG(I_{\dimp})))\mathrm{LB}(\lap(\UV))}}.$$ Finally, we can assume \(\lap(\UV)\geq1/2\) and that \(\sup_{\thv\in\UV}\|\nabla^3\lgd(\thv)\|_{\DFL}\sup_{\thv\in\UV}\|\DFL(\thv-\map)\|\leq 1/2\) for \(\DFL=I_{\dimp}/\normG(I_{\dimp})\). Thus we can get rid of the denominator. Note that the second of these assumptions is the same as we have made in~\eqref{TVsimpsimpcond}.

\section{Proofs for inverse problem}\label{app:inv}

%

\def\parm{q}
\subsection{Reduction to diagonal matrices}
\begin{proof}[Proof of Lemma~\ref{lma:D:prelim}]
A standard Riemann sum error bound gives
\begin{EQA}
|\Psi_{k\ell}-n\delta_{k\ell}|&=\left|\sum_{j=1}^n\psi_k(j/n)\psi_\ell(j/n)- n\int_0^1\psi_k(x)\psi_\ell(x)dx\right|\leq \|(\psi_k\psi_\ell)'\|_\infty\\
&\leq \|\psi_k\|_\infty\|\psi_\ell'\|_\infty+\|\psi_k'\|_\infty\|\psi_\ell\|_\infty\les k^\delta+\ell^\delta.
\end{EQA}To get the last inequality we used Assumption~\ref{assume:bdd} and that \(\sup_kk^{-\delta}\|\psi_k'\|_\infty\les1\). For a given \(\uv\in\R^{\dimp}\), we now have
\begin{EQA}
\bigg|\uv^\T \Psi\uv& -&n\|\uv\|^2\bigg|\les \sum_{k,\ell=1}^\dimp|u_k||u_\ell|( k^{\delta}+\ell^{\delta})= 2\sum_{\ell=1}^\dimp|u_\ell|\sum_{k=1}^\dimp|u_k|k^\delta\\
&\les&\left(\sum_{k=1}^\dimp|u_k|k^\delta\right)^2\leq \sum_{k=1}^{\dimp}|u_k|^2k^{\lambda}\sum_{k=1}^{\dimp}k^{2\delta-\lambda}\\
&\leq& (1+(\lambda-2\delta-1)^{-1})\uv^\T\diag(k^{\lambda})\uv.
\end{EQA} Here, we used~\eqref{ab-alpha-1}. Using that \(\lambda\) is bounded away from \(2\delta+1\) gives the desired result.
\end{proof}
\begin{proof}[Proof of Lemma~\ref{lma:hess}]
Recall from~\eqref{LL-radon} the form of \(\nabla^2L\), and recall that \(R_j^\T\theta = (\calR q_{\theta})(j/n)\). Since \(\sum_jR_jR_j^\T=R^\T R\) (because \(R_j\) are the rows of \(R\)), we have 
\begin{EQA}
&&c_hR^\T R\preceq \nabla^2L(\map)\preceq C_hR^\T R
\end{EQA}
Here, \(c_h=\inf\{ h''(t)\; : |t|\leq \|\calR q_{\map}\|_\infty\}\) and \(C_h=\sup\{ h''(t)\; : |t|\leq \|\calR q_{\map}\|_\infty\}\). Note that \(c_h>0\) because \(h\) is strictly convex. Also, thanks to Assumption~\ref{assume:Rqinfty}, we can assume \(c_h\) is bounded away from zero as \(n\to\infty\), for suitable \(\dimp\). Similarly, \(C_h\) is bounded above by a constant as \(n\to\infty\), for suitable \(\dimp\). We conclude that \(\DFL^2(\gamma_0)\asymp R^\T R+G_0^2\), proving~\eqref{ch-Ch-prelim}.

Let \(\Lambda=\diag(\lambda_k)\) and recall that \(\lambda_k\asymp k^{-2\beta}\). Also, recall from~\eqref{phi-psi} that \(\calR\phi_k=\psi_k/\sqrt{\lambda_k}\). Now, 
\begin{EQA}
(R^\T R)_{k\ell} &=& \sum_{j=1}^nR_{jk}R_{j \ell} = \sum_{j=1}^n\calR\phi_k(j/n)\calR\phi_\ell(j/n)\\
& = &\sqrt{\lambda_k}\sqrt{\lambda_\ell}\sum_{j=1}^n\psi_k(j/n)\psi_\ell(j/n) \\
&=& \sqrt{\lambda_k}\sqrt{\lambda_\ell}\Psi_{k\ell} = (\Lambda^{1/2}\Psi\Lambda^{1/2})_{k\ell},
\end{EQA} where \(\Psi\) is as in Lemma~\ref{lma:D:prelim}. Therefore, we have
\(R^\T R -n\Lambda = \Lambda^{1/2}(\Psi-nI_{\dimp})\Lambda^{1/2}\). Using the bounds on \(\Psi-nI_{\dimp}\) from Assumption~\ref{assume:ortho}, we conclude that
\begin{EQA}\label{Lkl}
&&-\Lambda^{1/2}\diag(k^{\lambda})\Lambda^{1/2}\precsim R^\T R-n\Lambda \precsim \Lambda^{1/2}\diag(k^{\lambda})\Lambda^{1/2}.
\end{EQA}
But note that
\begin{EQA}\label{Lkl2}
&&\Lambda^{1/2}\diag(k^{\lambda})\Lambda^{1/2} \asymp \diag(k^{\lambda-2\beta}) \preceq\diag(k^{2\gamma_0})=G_0^2,
\end{EQA}
since we assumed \(2\gamma_0\geq\lambda-2\beta\). Thus~\eqref{Lkl} and~\eqref{Lkl2} give \(-G_0^2 \precsim R^\T R -n\Lambda \precsim  G_0^2\), which implies \(R^\T R+G_0^2\asymp n\Lambda + G_0^2\).  Noting that \(\Lambda\asymp\diag(k^{-2\beta})\) concludes the proof.
\end{proof}

\subsection{Third derivative calculations}\label{app:sub:third}We start with some preliminary calculations which are independent of the particular choice of \(\DFL\). For a symmetric tensor  \(T\), the definition~\eqref{TDdef} of \(\|T\|_{\DFL}\) is equivalent to \(\sup_{\|\DFL\vv\|=1}\langle T, \vv^{\otimes3}\rangle\); see~\cite{symmtens}. Now, fix some \(\thv\) and note that 
\begin{EQA}
&&\nabla^3\lgd(\thv)=\sum_{j=1}^nh'''(R_j^\T\thv)R_j^{\otimes3}.
\end{EQA} Recall also that \(R_j^\T\thv=(\calR q_{\thv})(j/n)\). Therefore
\begin{EQA}\label{nabla3Lthv}
\|\nabla^3\lgd(\thv)\|_{\DFL}&=& \sup_{\|\DFL\vv\|=1}\left|\sum_{j=1}^nh'''((\calR q_{\thv})(j/n))(R_j^\T\vv)^3\right|\\
&\leq &\|h'''(\calR q_{\thv})\|_\infty\sup_{\|\DFL\vv\|=1}\sum_{j=1}^n|R_j^\T\vv|^3.
\end{EQA} 
Thus for a set \(\UV\), we have
\begin{EQA}\label{nabla3Lthv1}
&&\sup_{\thv\in\map+\UV}\|\nabla^3\lgd(\thv)\|_{\DFL}\leq \sup_{\thv\in\map+\UV}\|h'''(\calR q_{\thv})\|_\infty\sup_{\|\DFL\vv\|=1}\sum_{j=1}^n|R_j^\T\vv|^3,
\end{EQA} and if \(\|h'''\|_\infty<\infty\), then
\begin{EQA}\label{htriple}
&&\sup_{\thv\in\R^{\dimp}}\|\nabla^3\lgd(\thv)\|_{\DFL}\leq \|h'''\|_\infty\sup_{\|\DFL\vv\|=1}\sum_{j=1}^n|R_j^\T\vv|^3.
\end{EQA}
Next, we observe that
\begin{EQA}\label{R1R2bd}
\sup_{\|\DFL\vv\|=1}\sum_{j=1}^n|R_j^\T\vv|^3 &\leq& \sup_{\|\DFL\vv\|=1,j\in[n]}|R_j^\T\vv|\sup_{\|\DFL\vv\|=1}\vv^\T\left(\sum_{j=1}^nR_jR_j^\T\right)\vv\\
&\leq &\sup_{\|\DFL\vv\|=1}\|\calR q_{\vv}\|_\infty \sup_{\|\DFL\vv\|=1}\vv^\T R^\T R\vv.
\end{EQA} We also record a coarser bound, which will turn out to be tight when \(\DFL=I_{\dimp}/\normG(I_{\dimp})\):
\begin{EQA}\label{R1onlybd}
&&\sup_{\|\DFL\vv\|=1}\sum_{j=1}^n|R_j^\T\vv|^3  \leq n\sup_{\|\DFL\vv\|=1}\|\calR q_{\vv}\|^3_\infty.
\end{EQA}
We now specialize to the relevant choices of \(\DFL\).
\begin{lemma}For all \(2\gamma_0+2\beta\geq\lambda\), we have
\begin{EQA}\label{R1R2Dgam}
&&\sup_{\|\DFL(\gamma_0)\vv\|=1}\|\calR q_{\vv}\|_\infty\les S_{\dltwu}(\gamma_0),\qquad \sup_{\|\DFL(\gamma_0)\vv\|=1}\vv^\T R^\T R\vv\les 1.
\end{EQA} Furthermore, for all \(\beta\neq1/2\), 
\begin{EQA}\label{R1Did}
&&\sup_{\|\vv\|=1}\|\calR q_{\vv}\|_\infty\les\dimp^{(1-2\beta)\vee0}.
\end{EQA}
\end{lemma}
\begin{proof}Using \(\calR\phi_k = \sqrt{\lambda_k}\psi_k\) by~\eqref{phi-psi}, the fact that \(\lambda_k\asymp k^{-\beta}\) by~\eqref{lambda2beta}, and that \(\|\psi_k\|_\infty \les 1\) by Assumption~\ref{assume:bdd}, we have
\begin{EQA}\label{bdonT}
\|\calR q_{\vv}\|_\infty&\leq& \sum_{k=1}^{\dimp}|v_k|\|\calR \phi_k\|_\infty\les \sum_{k=1}^{\dimp}|v_k|k^{-\beta}\\
&\leq &\sqrt{\sum_{k=1}^{\dimp}\frac{k^{-2\beta}}{n/k^{2\beta}+k^{2\gamma_0}}}\sqrt{\sum_{k=1}^{\dimp}v_k^2\left(\frac{n}{k^{2\beta}}+k^{2\gamma_0}\right)}\\
&\les&S_{\dltwu}(\gamma_0)\|\DFL(\gamma_0)\vv\|.
\end{EQA} To get the last line, we used that \(\DFL(\gamma_0)^2\succsim\diag(n/k^{2\beta}+k^{2\gamma_0})\) by Lemma~\ref{lma:hess}. Taking the supremum over \(\|\DFL(\gamma_0)\vv\|=1\) concludes the proof of the first inequality in~\eqref{R1R2Dgam}. For the second inequality, we use that \(R^\T R\prec R^\T R+G_0^2\preceq\DFL(\gamma_0)^2\), again by Lemma~\ref{lma:hess}.

For~\eqref{R1Did}, we have, similarly to~\eqref{bdonT}, that
\begin{EQA}
\sup_{\|\vv\|=1}\|\calR q_{\vv}\|_\infty&\les &\sup_{\|\vv\|=1}\sum_{k=1}^{\dimp}|v_k|k^{-\beta}\leq \sqrt{\sum_{k=1}^{\dimp}k^{-2\beta}} \asymp \dimp^{(1-2\beta)\vee0}
\end{EQA} for all \(\beta\neq1/2\). This is by~\eqref{asympsum1} and~\eqref{ab-alpha-1}.
\end{proof}

We conclude from~\eqref{R1R2bd} and~\eqref{R1R2Dgam} that for all \(\gamma_0\) such that \(2\beta+2\gamma_0\geq\lambda\), it holds
\begin{EQA}\label{R3-final-Dgam}
&&\sup_{\|\DFL(\gamma_0)\vv\|=1}\sum_{j=1}^n|R_j^\T\vv|^3 \les S_{\dltwu}(\gamma_0).
\end{EQA} If \(\beta>1/2\), then~\eqref{R1onlybd} and~\eqref{R1Did} give
\begin{EQA}\label{R3-final-Did}
&&\sup_{\|\vv\|=1}\sum_{j=1}^n|R_j^\T\vv|^3 \les n.
\end{EQA}

\subsubsection{Tightness of the operator norm bounds}\label{app:tight}
When \(\|h'''\|_\infty<\infty\), the first step in our operator norm bounds is the initial bound~\eqref{htriple}, which follows from passing the absolute value in the first line of~\eqref{nabla3Lthv} into the sum and then bounding \(|h'''(\cdot)|\) by its \(L^\infty\) norm. The next and only remaining step is to bound  \(\sup_{\|\DFL\vv\|=1}\sum_{j=1}^n|R_j^\T\vv|^3\), as we do in~\eqref{R3-final-Dgam} and~\eqref{R3-final-Did}. We now show that these two bounds are tight, under the simplifying assumption that the \(\psi_k\) are an exact Fourier basis. 

Namely, assume the \(\psi_k\) are given by the orthonormal functions \(\psi_k(x)=\sqrt2\cos(\pi kx)\), \(k=1,2,3,\dots\). Thus \(\calR\) is determined by the equations \(\calR \phi_k=\sqrt{\lambda_k}\psi_k\) for some orthonormal basis of functions \(\phi_k\), \(k=1,2,3,\dots\) in \(L^2[0,1]\).
\begin{proposition}\label{prop:tighttau3}Let \(\calR\) be as described above, with \(\lambda_k\asymp k^{-2\beta}\), \(k=1,2,3,\dots\). Suppose \(2\gamma_0+2\beta\geq\lambda\) and \(2\gamma_0+2\beta-1\) is bounded below by a positive constant depending only on \(\gamma\) and \(\beta\). Then 
\begin{EQA}\label{eq:tighttau3}
&&\sup_{\|\DFL(\gamma_0)\vv\|=1}\sum_{j=1}^n|R_j^\T v|^3 \asymp S_{\dltwu}(\gamma_0) \asymp \sqrt{\frac{m_0(\gamma_0)\wedge\dimp}{n}},
\end{EQA} where \(m_0(\gamma_0)=n^{1/(2\beta+2\gamma_0)}\). For the same \(\calR\) and \(\lambda_k\)'s, let \(\DFL=I_{\dimp}/\normG(I_{\dimp})\). If \(\beta>1/2\), then
\begin{EQA}\label{DFLI}
&&\sup_{\|\DFL\vv\|=1}\sum_{j=1}^n|R_j^\T \vv|^3 =\normG(I_{\dimp})^3\sup_{\|\vv\|=1}\sum_{j=1}^n|R_j^\T\vv|^3\asymp \normG(I_{\dimp})^3 n.
\end{EQA}
\end{proposition}Here, recall that the constants of proportionality may depend on \(\gamma\) and \(\beta\). The value \(\gamma_0=\gamma\) satisfies the assumptions of this proposition, since \(2\gamma+2\beta-1>0\). Recalling that \(\DFL(\gamma)=\DFLG\) and noting that \(m_0(\gamma)=m\), we conclude from~\eqref{eq:tighttau3} that the bound on \(\|\nabla^3\lgd\|_{\DFLG,\infty}\) in Lemma~\ref{lma:compare} is tight up to the initial step. Next, we show that \(\gamma_0=\gamma_0^*=\gamma-1/2-1/(2m)\) also satisfies the assumption that \(2\gamma_0+2\beta-1\) is bounded below by a positive constant depending only on \(\gamma\) and \(\beta\), provided \(n\) is large enough. Indeed, suppose \(n>(\beta+\gamma-1)^{-2\beta-2\gamma}\), in which case \(m^{-1}\leq\beta+\gamma-1\). Then \(2\gamma_0+2\beta-1= 2\gamma+2\beta-2 -1/m\geq\gamma+\beta-1\), which is a positive constant. Thus we see the stated bound on \(\|\nabla^3\lgd\|_{\DFL(\gamma_0^*),\infty}\) in Lemma~\ref{lma:compare} is tight up to the initial step described above. Finally,~\eqref{DFLI} can be combined with the fact that
\begin{EQA}
&&\normG(I_{\dimp})\asymp (m\wedge\dimp)^\beta/\sqrt n,
\end{EQA}
shown in~\eqref{normG-id}, to conclude that the bound on \(\|\nabla^3\lgd\|_{I_{\dimp}/\normG(I_{\dimp}),\infty}\) in Lemma~\ref{lma:compare} is tight up to the initial step. 
\begin{proof} First note that
$$\DFL_1\asymp\DFL_2\quad\implies\quad\sup_{\|\DFL_1\vv\|=1}\sum_{j=1}^n|R_j^\T \vv|^3\asymp \sup_{\|\DFL_2\vv\|=1}\sum_{j=1}^n|R_j^\T \vv|^3.$$ Thus we can assume without loss of generality that \(\DFL(\gamma_0)^2= \diag(n/k^{2\beta}+k^{2\gamma_0})\), since the condition \(2\gamma_0+2\beta\geq\lambda\) from Lemma~\ref{lma:hess} is satisfied. Second, recall that \(R_j^\T\vv=(\calR q_{\vv})(j/n)=\sum_{k=1}^{\dimp}v_k\sqrt{\lambda_k}\psi_k(j/n) \asymp \sum_{k=1}^{\dimp}v_kk^{-\beta}\psi_k(j/n)\), where the constant of proportionality does not depend on \(j\).  Thus we can assume \(\lambda_k=k^{-\beta}\).

Let \(\bar m_0=\min(m_0(\gamma_0),\dimp)\), where \(m_0(\gamma_0)=n^{1/(2\beta+2\gamma_0)}\). Note \(\bar m_0\geq1\). Combining Lemmas~\ref{lma:bet-gam-firstcase} and~\ref{lma:bet-gam-lastcase}, and the assumption that \(2\gamma_0+2\beta-1\) is bounded below away from zero (by a \(\beta,\gamma\)-dependent constant), we have that \(S_{\dltwu}^2(\gamma_0)\asymp\bar m_0/n\).

Furthermore,~\eqref{R3-final-Dgam} shows that \(\sup_{\|\DFL(\gamma_0)\vv\|=1}\sum_{j=1}^n|(\calR q_{\vv})(j/n)|^3\les S_{\dltwu}(\gamma_0)\). Therefore, to finish the proof of~\eqref{eq:tighttau3}, it remains to exhibit \(\vv\) such that \(\|\DFL(\gamma_0)\vv\|=1\) and \(\sum_{j=1}^n|(\calR q_{\vv})(j/n)|^3\gtrsim \sqrt{\bar m_0/n}\). Let \(v_k = Ck^{\beta}/\sqrt{\bar m_0n}\) for all \(k=1,\dots, \bar m_0\) and \(v_k=0\) for all \(k>\bar m_0\). We choose \(C\) so that \(\|\DFL(\gamma_0)\vv\|=1\). Specifically, we have
\begin{EQA}
&&1=\|\DFL(\gamma_0)\vv\|^2 = \frac{C^2}{\bar m_0n}\sum_{k=1}^{\bar m_0}(n+k^{2\gamma_0+2\beta}) \leq 2C^2,
\end{EQA} using that \(k^{2\gamma_0+2\beta}\leq m_0^{2\gamma_0+2\beta}=n\) for all \(k\leq\bar m_0\). Thus \(C\geq1/\sqrt2\). Now, we have
\begin{EQA}
&&\calR q_{\vv}(x) = \sum_{k=1}^{\bar m_0}v_k\calR\phi_k(x) = \sum_{k=1}^{\bar m_0}v_kk^{-\beta}\psi_k(x) = \frac{C\sqrt2}{\sqrt{\bar m_0n}} \sum_{k=1}^{\bar m_0}\cos(\pi kx).
\end{EQA} Now, fix \(x\in [0,(4\pi(\bar m_0+1))^{-1}]\). Using the inequality \(|\cos a-1|\leq |a|\) for all \(a\), we have
\begin{EQA}
\left|\left( \sum_{k=1}^{\bar m_0}\cos(\pi kx)\right)-\bar m_0\right|&\leq&\sum_{k=1}^{\bar m_0}|\cos(\pi kx)-1|\leq \sum_{k=1}^{\bar m_0}\pi kx\\
&\leq &\pi(\bar m_0+1)^2x \leq (\bar m_0+1)/4\leq\bar m_0/2.
\end{EQA}
We conclude that 
\begin{EQA}
&&|\calR q_{\vv}(x)| \geq \frac{C\sqrt2}{\sqrt{\bar m_0n}} \frac{\bar m_0}{2} \geq \frac 12\sqrt{\frac{\bar m_0}{n}}\quad\forall x\in [0,(4\pi \bar m_0)^{-1}].
\end{EQA} From here, we get that
\begin{EQA}
&&\sum_{j=1}^n|(\calR q_{\vv})(j/n)|^3\geq \sum_{j=1}^{ n/(4\pi \bar m_0)}|\calR q_{\vv}(j/n)|^3\gtrsim \frac{n}{\bar m_0}\frac{\bar m_0\sqrt{\bar m_0}}{n\sqrt n} =\sqrt{\frac{\bar m_0}{n}}.
\end{EQA} This concludes the proof of the first statement.

Next, suppose \(\DFL=I_{\dimp}/\normG(I_{\dimp})\). We have from~\eqref{R3-final-Did} that \(\sup_{\|\vv\|=1}\sum_{j=1}^n|(\calR q_{\vv})(j/n)|^3\les n\) when \(\beta\) is bounded below away from \(1/2\). It remains to show this upper bound is tight. Take \(\vv=(1,0,\dots,0)\), so that \(\calR q_{\vv}(x)=\sqrt2\cos(\pi x)\). Then there are \(\mathcal O(n)\) values of \(j\) for which \((\calR q_{\vv})(j/n)\geq1\), by the same arguments as above. This shows \(\sup_{\|\vv\|=1}\sum_{j=1}^n|(\calR q_{\vv})(j/n)|^3\geq cn\), as desired.
\end{proof}



\subsection{Proof of Lemma~\ref{lma:effdim:radon}, Theorem~\ref{corr:TV:radon}, Corollary~\ref{corr:TVgam0star}, Lemma~\ref{lma:compare}}\label{app:sec3:proofs}
\begin{proof}[Proof of Lemma~\ref{lma:effdim:radon}]
Since \(\normG(\DFL(\gamma_0))=1\), we have \(\dimL(\DFL(\gamma_0)) =\Tr(\DFLG^{-1}\DFL^2(\gamma_0)\DFLG^{-1})\). Note that \(2\gamma_0+2\beta\geq\lambda\) implies also \(2\gamma+2\beta\geq\lambda\), and Lemma~\ref{lma:hess} then gives that \(\DFL^2(\gamma_0)\asymp\diag(n/k^{2\beta}+k^{2\gamma_0})\) and \(\DFL^2(\gamma)=\DFLG^2\asymp\diag(n/k^{2\beta}+k^{2\gamma})\). Finally, we apply Lemma~\ref{lma:AB} to conclude  \(\dimL(\DFL(\gamma_0))\asymp S_{\dimL}(\gamma_0)\). For the rest of the proof, we use \(\DFL\) and \(S_{\dltwu}\) as shorthand for \(\DFL(\gamma_0)\) and \(S_{\dltwu}(\gamma_0)\), respectively. Combining~\eqref{nabla3Lthv1} and~\eqref{R3-final-Dgam}, we have
\begin{EQA}
&&\sup_{\|\DFL\uv\|\leq\rrdc}\|\nabla^3\lgd(\map+\uv)\|_{\DFL}\leq \sup_{\|\DFL\uv\|\leq\rrdc}\|h'''(\calR q_{\map+\uv})\|_\infty S_{\dltwu}.
\end{EQA}
Using that \(\calR q_{\map+\uv}=\calR q_{\map} + \calR q_{\uv}\), we have
\begin{EQA}\label{bd1t3}
\sup_{\|\DFL\uv\|\leq\rrdc}\|h'''(\calR q_{\map+\uv})\|_\infty  &\leq &\sup\{|h'''(t)|\,:\,|t|\leq\|\calR q_{\map}\|_\infty + \sup_{\|\DFL\uv\|\leq\rrdc}\|\calR q_{\uv}\|_\infty\}\\
&\leq &\sup\{|h'''(t)|\,:\,|t|\leq C+C\rrdc S_{\dltwu}\} \\
&\leq & D_3(C(1\vee \rrdL S_{\dltwu})).
\end{EQA} Here, we used boundedness of \(\|\calR q_{\map}\|_\infty\) by Assumption~\ref{assume:Rqinfty}, as well as the first inequality in~\eqref{R1R2Dgam}. Combining the above two displays concludes the proof.
\end{proof}
\begin{proof}[Proof of Theorem~\ref{corr:TV:radon}]We apply Corollary~\ref{corr:DD0}, with \(\DFL=\DFL(\gamma_0)\). Recall that (1) and (2) are satisfied by construction. The condition (3) is satisfied if\\ \(\sup_{\|\DFL(\gamma_0)\uv\|\leq \rrdc}\|\nabla^3\lgd(\map+\uv)\|_{\DFL(\gamma_0)}\rrdc \leq 1/2\) for some \(\rrdc\geq 6\sqrt{\dimL(\DFL(\gamma_0))}\geq3\sqrt{\dimL(\DFL(\gamma_0))}+3\). The bounds from Lemma~\ref{lma:effdim:radon} show that this is satisfied as long as
\begin{EQA}\label{eq:r1}
&&\rrdL \geq C_1\sqrt{S_{\dimL}},\\
&&D_3(C(1\vee \rrdL S_{\dltwu}))\rrdL S_{\dltwu} \leq c_2\label{eq:r2}
\end{EQA} for some \(C_1\) large enough and \(c_2\) small enough. For any \(\rrdL\) satisfying these conditions,~\eqref{eq:TVsimpsimp} and Lemma~\eqref{lma:effdim:radon} give
\begin{EQA}\label{TV-radon-1}
&&\TV( \PfL,\,\lap)\les D_3(C(1\vee \rrdL S_{\dltwu}))S_{\dltwu}S_{\dimL}+ \ex^{-\rrdL^2/12}.
\end{EQA}
Next, note that if \(S_{\dltwu}S_{\dimL}\) is small enough then so is \(S_{\dltwu}\), since \(S_{\dimL}\geq1\). Thus we can assume in particular that \(S_{\dltwu}\leq1\). Take \(\rrdL=1/\sqrt{S_{\dltwu}}\), which gives \(D_3(C(1\vee \rrdL S_{\dltwu}))=D_3(C(1\vee \sqrt{S_{\dltwu}}))=D_3(C)\), which is a constant, since \(h\) is fixed. For this \(\rrdL\), the condition~\eqref{eq:r1} reduces to \(S_{\dltwu}S_{\dimL}\) being small enough. The condition~\eqref{eq:r2} reduces to \(S_{\dltwu}\) being small enough, which is weaker. We conclude the conditions are satisfied. Plugging in \(\rrdL=1/\sqrt{S_{\dltwu}}\) to~\eqref{TV-radon-1} gives 
\begin{EQA}\label{app:TV-radon-2}
&&\TV( \PfL,\,\lap)\les D_3(C)S_{\dltwu}S_{\dimL}+ \ex^{-1/(12S_{\dltwu})}\les  S_{\dltwu}S_{\dimL}.
\end{EQA}
\end{proof}
\begin{remark}
The second bound in~\eqref{app:TV-radon-2} is designed for the case that \(h'''\) is nonnegligible. However if e.g. \(h'''=0\), we obtain from the first bound that \(\TV( \PfL,\lap)\les \exp(-1/(12S_{\dltwu}))\), which is much tighter. Note that when \(h\) is quadratic, \( \PfL\) is Gaussian, so the LA is exact and the TV distance is zero. 
\end{remark}
\begin{proof}[Proof of Corollary~\ref{corr:TVgam0star}]We apply Theorem~\ref{corr:TV:radon} with \(\gamma_0=\gamma_0^*=\gamma-1/2-1/(2m)\), which requires that \(2\gamma_0^*+2\beta=2\gamma+2\beta-1-1/m\geq\lambda\) is satisfied. This is true for \(n\) large enough, since we have that \(2\gamma+2\beta-1-\lambda>0\).We conclude that
\begin{EQA}
&&\TV(\PfL,\lap)\les S_{\dltwu}(\gamma_0^*)S_{\dimL}(\gamma_0^*),
\end{EQA} whenever this quantity is small enough. To conclude, we apply Lemma~\ref{lma:Sranges}, whose assumptions are also satisfied for \(n\) large enough.
\end{proof}
\begin{proof}[Proof of Lemma~\ref{lma:compare}]
Lemma~\ref{lma:effdim:radon} gives \(\dimL(\DFL(\gamma_0^*))\asymp S_{\dimL}(\gamma_0^*)\), and Lemma~\ref{lma:Sranges} gives \(S_{\dimL}(\gamma_0^*)\asymp m\wedge\dimp\). The fact that \(\dimL(\DFLG^2)=\dimp\) follows from the definition~\eqref{eq:effdim} of \(\dimL\). Combining~\eqref{htriple} and~\eqref{R3-final-Dgam}, we have \(\sup_{\thv\in\R^{\dimp}}\|\nabla^3\lgd(\thv)\|_{\DFL(\gamma_0)}\les \|h'''\|_\infty S_{\dltwu}(\gamma_0)\les S_{\dltwu}(\gamma_0)\). When \(\gamma_0=\gamma_0^*\), we conclude by Lemma~\ref{lma:Sranges}. When \(\gamma_0=\gamma\), we obtain \(\sup_{\thv\in\R^{\dimp}}\|\nabla^3\lgd(\thv)\|_{\DFLG}\les S_{\dltwu}(\gamma)\les \sqrt{(m\wedge\dimp)/n}\). The last inequality follows by combining Lemmas~\ref{lma:bet-gam-firstcase} and~\ref{lma:bet-gam-lastcase}, and recalling that \(m_0=m\) when \(\gamma_0=\gamma\).

We now turn to the analysis of the case \(\DFL=I_{\dimp}/\normG(I_{\dimp})\). First, using that \(\DFLG^2=\DFL^2(\gamma)\asymp\diag(n/k^{2\beta}+k^{2\gamma})\) by Lemma~\ref{lma:hess},  have 
\begin{EQA}\label{normG-id}
&&\normG(I_{\dimp})^2=\|\DFLG^{-2}\| \asymp \|\diag((n/k^{2\beta}+k^{2\gamma})^{-1})\|\asymp (m\wedge\dimp)^{2\beta}/n.
\end{EQA} Next,
\begin{EQA}
\Tr(\DFLG^{-2})& \asymp&\Tr(\diag(n/k^{2\beta}+k^{2\gamma})^{-1})=\sum_{k=1}^{\dimp}\frac{k^{2\beta}}{n+k^{2\gamma+2\beta}}\\
&\asymp & \sum_{k=1}^{m\wedge\dimp}k^{2\beta}/n + \sum_{k=1+(m\wedge\dimp)}^{\dimp}k^{-2\gamma}\\
&\asymp& (m\wedge\dimp)^{2\beta+1}/n + \sum_{k=1+(m\wedge\dimp)}^{\dimp}k^{-2\gamma}
\end{EQA} Here, we have used~\eqref{splitsum} and~\eqref{asympsum1}. Next, we use that \(\gamma>1/2\) and the upper bound in~\eqref{tail-sum}. This gives
\begin{EQA}
&&(m\wedge\dimp)^{2\beta+1}/n  \les \Tr(\DFLG^{-2})\les (m\wedge\dimp)^{2\beta+1}/n + m^{1-2\gamma}\Ind\{m<\dimp\}.
\end{EQA}
Using the definition of \(m\), this implies  \(\Tr(\DFLG^{-2})\asymp (m\wedge\dimp)^{2\beta+1}/n\). Thus
$$\dimL(I_{\dimp}/\normG(I_{\dimp}))= \Tr(\DFLG^{-2})/\|\DFLG^{-2}\|\asymp m\wedge\dimp.$$ Finally,~\eqref{htriple} and~\eqref{R3-final-Did} give \(\sup_{\thv\in\R^{\dimp}}\|\nabla^3\lgd(\thv)\| \les n\), and therefore 
\begin{EQA}
\sup_{\thv\in\R^{\dimp}}\|\nabla^3\lgd(\thv)\|_{I_{\dimp}/\normG(I_{\dimp})}&\les & \normG(I_{\dimp})^3n \asymp n^{-1/2}(m\wedge\dimp)^{3\beta}.
\end{EQA}
\end{proof}


\subsection{Bounds on \(S_{\dimL}\) and \(S_{\dltwu}\), and proof of Lemma~\ref{lma:Sranges}}\label{app:Sbounds}
Recall that 
\begin{EQA}\label{app:Sdefs}
&&S_{\dimL}=\sum_{k=1}^{\dimp}\frac{n+k^{2\gamma_0+2\beta}}{n+k^{2\gamma+2\beta}},\qquad S_{\dltwu}^2= \sum_{k=1}^{\dimp}\frac{1}{n+k^{2\beta+2\gamma_0}}.
\end{EQA}
Furthermore, recall the definition of \(m,m_0\) from~\eqref{def:m-m0}.
\begin{lemma}\label{lma:bet-gam-firstcase} For all \(\gamma_0\leq\gamma\) we have
\begin{EQA}\label{app:ub-lb-univ}
\frac12(m\wedge\dimp)&\leq& S_{\dimL}\leq\dimp,\qquad \frac{m_0\wedge\dimp}{2n}\leq S_{\dltwu}^2\leq\frac{\dimp}{n}.
\end{EQA}
\end{lemma}
\begin{proof}The upper bound on \(S_{\dimL}\) is immediate from the definition and the fact that \(\gamma_0\leq\gamma\). The lower bound stems from the fact that \(n+k^{2\beta+2\gamma}\leq 2n\) and \(n+k^{2\beta+2\gamma_0}\geq n\) for all \(k=1,\dots,m\wedge\dimp\). The upper bound on \(S_{\dltwu}\) uses that \(1/(n+k^{2\beta+2\gamma_0})\leq 1/n\) for all \(k=1,\dots,\dimp\). For the lower bound, we use that \(n+k^{2\beta+2\gamma_0}\leq 2n\) for all \(k=1,\dots,m_0\wedge\dimp\).
\end{proof}
\begin{lemma}\label{lma:bet-gam-dimL}Suppose \(2\beta+2\gamma>1\) and \(\dimp>m\). If \(\gamma_0<\gamma-1/2\), then
\begin{EQA}
&&S_{\dimL}\leq \left(1+(2\gamma+2\beta-1)^{-1}+(2\gamma-2\gamma_0-1)^{-1}m^{2\gamma_0-2\gamma}\right)m.
\end{EQA}
\end{lemma} 
\begin{proof}We have
\begin{EQA}
S_{\dimL}&\leq &m + \sum_{k=m+1}^{\dimp}\frac{n+k^{2\gamma_0+2\beta}}{k^{2\gamma+2\beta}}\\
&=&m + n\sum_{k=m+1}^{\dimp}k^{-2\beta-2\gamma}+\sum_{k=m+1}^{\dimp}k^{2\gamma_0-2\gamma}\\
&\leq &m + \frac{n}{2\beta+2\gamma-1}m^{1-2\beta-2\gamma} + \frac{1}{2\gamma-2\gamma_0-1}m^{1+2\gamma_0-2\gamma}\\
&= & \left(1+(2\gamma+2\beta-1)^{-1}+(2\gamma-2\gamma_0-1)^{-1}m^{2\gamma_0-2\gamma}\right)m.
\end{EQA}
To get the third line we used~\eqref{tail-sum}. To get the fourth line we used \(nm^{1-2\beta-2\gamma}=m\). 
\end{proof}

Next, we analyze \(S_{\dltwu}\). 
\begin{lemma}\label{lma:bet-gam-lastcase}If \(\dimp>m_0\) and \(\gamma_0>1/2-\beta\), then
\begin{EQA}\label{ub-Stau}
&& S_{\dltwu}^2\leq \left(1+(2\beta+2\gamma_0-1)^{-1}\right)\frac{m_0}{n}.
\end{EQA}
\end{lemma}
\begin{proof}We use the inequality \(n+k^{2\beta+2\gamma_0}\geq n\) for all \(k=1,\dots,m_0\) and the inequality \(n+k^{2\beta+2\gamma_0}\geq k^{2\beta+2\gamma_0}\) for all \(k=m_0+1,\dots,\dimp\). Using this and~\eqref{tail-sum}, we have
\begin{EQA}
S_{\dltwu}^2&\leq &\frac{m_0}{n} + \sum_{k>m_0}k^{-2\beta-2\gamma_0} \leq \frac{m_0}{n}  + (2\beta+2\gamma_0-1)^{-1}m_0^{1-2\beta-2\gamma_0}\\ &=&\left(1+(2\beta+2\gamma_0-1)^{-1}\right)\frac{m_0}{n}.
\end{EQA} \end{proof}
Combining the above three lemmas gives that when \(2\beta+2\gamma>1\) and \(\gamma_0<\gamma-1/2\), we have
\begin{EQA}\label{SdimLrange-2}
&& \frac{m\wedge\dimp}{2}\leq S_{\dimL}\leq \big(1 + (2\beta+2\gamma-1)^{-1}+(2\gamma-2\gamma_0-1)^{-1}m^{2\gamma_0-2\gamma}\big)(m\wedge\dimp).\end{EQA} When \(\gamma_0>\frac12-\beta\), we have  
 \begin{EQA}
 &&\frac{m_0\wedge\dimp}{2n}\leq S_{\dltwu}^2\leq  \big(1+(2\beta+2\gamma_0-1)^{-1}\big)\frac{m_0\wedge\dimp}{n}.\label{Stau3range-2}
\end{EQA} 
We now use~\eqref{SdimLrange-2} and~\eqref{Stau3range-2} to choose \(\gamma_0\). Note that the bound on \(S_{\dimp}\) increases with \(\gamma_0\), but the dependence is mild: through the coefficient multiplying \(m\wedge\dimp\) and not through a power. Meanwhile, the bound on \(S_{\dltwu}\) decreases as \(\gamma_0\) increases, and the dependence is strong: through the power of \(n\) defining \(m_0\). Thus the strategy to choose \(\gamma_0\) is to take it as large as possible, without \(c=(2\gamma-2\gamma_0-1)^{-1}m^{2\gamma_0-2\gamma}\) blowing up. In particular, it is reasonable to choose \(\gamma_0\) to ensure \(c\) is bounded by a constant. We choose \(\gamma_0=\gamma_0^*=\gamma-1/2-1/2m\), which yields \(c\leq1\). 
This gives
\begin{EQA}
\frac{m\wedge\dimp}{2}&\leq &S_{\dimL}\leq (2 + (2\beta+2\gamma-1)^{-1})(m\wedge\dimp),\label{Sdimfinal}\\
\frac{m_0^*\wedge\dimp}{2n}&\leq &S_{\dltwu}^2\leq  \big(1+(2\beta+2\gamma-2-1/m)^{-1}\big)\frac{m_0^*\wedge\dimp}{n}.\label{Staufinal}
\end{EQA}
Finally, to finish the proof of Lemma~\ref{lma:Sranges}, we should have that \(2\beta+2\gamma-2-1/m>0\), and in fact, \(1/m\) should be small enough that this quantity can be lower bounded in an \(n\)-independent way, so that the constant of proportionality in the upper bound~\eqref{Staufinal} does not depend on \(n\). It suffices to assume \(2\beta+2\gamma>2\) and \(m^{-1}\leq\beta+\gamma-1\) or equivalently, that \(n>(\beta+\gamma-1)^{-2\beta-2\gamma}\).

The upper bound in~\eqref{Staufinal} then becomes
\begin{EQA}
&&S_{\dltwu}^2\leq  \big(1+(\beta+\gamma-1)^{-1}\big)\frac{m_0^*\wedge\dimp}{n}
\end{EQA}
This proves Lemma~\ref{lma:Sranges}. 

\section{Eigenvalues, eigenfunctions for Volterra-like operators}\label{app:eigen}
We consider integral operators \(\calR:L^2[0,1]\to L^2[0,1]\) defined by 
\begin{EQA}\label{calR-def}
&&g=\calR f \quad \Leftrightarrow\quad ag' + bg = f,\quad g(0)=0,
\end{EQA}  In other words, \(\calR f\) is defined as the unique function \(g\) which solves the initial value ODE on the right.  As noted in Proposition~\ref{thm:psik0}, \(\calR\) can equivalently be expressed as the Volterra-like integral operator~\eqref{calR-int-def0}.

In this section, we prove Proposition~\ref{thm:psik0}, restated here for convenience.
\begin{theorem}\label{thm:psik}
Let \(\calR\) be as in~\eqref{calR-def}, where \(a\in C^2[0,1]\), \(b\in C^1[0,1]\), and there is \(a_0>0\) such that \(a(x)>a_0>0\) for all \(x\in[0,1]\). Then \(\calR\) is compact. Let \(\psi_k\) be any complete orthonormal basis of \(L^2[0,1]\) of eigenfunctions of \(\calR\calR^\T\), and \(\lambda_k\) be the corresponding eigenvalues. Then \(\lambda_k\asymp k^{-2}\). Furthermore, the \(\psi_k\) satisfy
\begin{EQA}\label{psik-assume}
&&\|\psi_k\|_\infty \leq \Psi_0,\qquad \|\psi_k'\|_\infty \leq \Psi_1 k^\delta
\end{EQA} with \(\delta=1\). The constants \(\Psi_0,\Psi_1\) depend only on the functions \(a\) and \(b\).
\end{theorem}
Note that since \(a,b\) are continuous and \(a\) is bounded below away from zero, \(\calR\) is a Hilbert-Schmidt operator and therefore compact. To prove the rest of the theorem, we start by computing \(\calR\calR^\T\). 
\begin{lemma}\label{lma:RRT}Let \(\calR\) be as in~\eqref{calR-def}. Then the operator \(\calR\calR^\T\) is defined as follows:
\begin{EQA}\label{calRRT-def}
&&h=\calR\calR^\T f \quad \Leftrightarrow\quad \text{\(h\) solves}\quad
\begin{cases}-(ph')' + qh = f, \quad x\in(0,1),\\
h(0)=0,\quad a(1)h'(1)+b(1)h(1)=0,\end{cases}
\end{EQA} where \(p=a^2\) and \(q=b^2-(ab)'\).
\end{lemma} See the end of this subsection for the proof. Now, since \(\calR\calR^\T\) is nonnegative definite by construction, all of its eigenvalues \(\lambda_k\) must be nonnegative. Furthermore, if \(\calR\calR^\T f = 0=h\), then~\eqref{calRRT-def} implies \(f=0\). Thus \(\lambda_k>0\) for all \(k\). Eigenfunctions \(\psi_k\) of \(\calR\calR^\T\) are equivalently eigenfunctions of the Sturm-Liouville differential operator, with
\begin{EQA}
\begin{cases}-(p\psi_k')' + q\psi_k = \lambda_k^{-1}\psi_k \quad x\in(0,1),\\
\psi_k(0)=0,\quad a(1)\psi_k'(1)+b(1)\psi_k(1)=0,\end{cases} &&
\end{EQA}
We now use Theorem 4.3.1 of~\cite{zettl2005sturm}. Note that the above boundary conditions (BCs) are ``separated", i.e. in the form (4.2.1),(4.2.2) of~\cite{zettl2005sturm}, which is a special case of the BC required by Theorem 4.3.1, given in (4.1.3),(4.1.4). Part (7) of the theorem shows that \(\lambda_k^{-1}\to\infty\) at rate \(k^2\). Specifically, we have
\begin{EQA}\label{lambdak-asymp}
&&k^2\lambda_k\to\pi^{-2}\left(\int_0^1a(t)^{-1}dt\right)^{2},\qquad k\to\infty.
\end{EQA}
We now study the asymptotics of the eigenfunctions \(\psi_k\), in order to check the condition~\eqref{psik-assume}. To do so, we use the Liouville transformation; see Chapter 7 of~\cite{everitt2005sturm}. Let \(T=\int_0^1a(x)^{-1}dx\). Define the transformation \(t(x)=\int_0^xa(s)^{-1}dx\), so that \(t\) increases monotonically from \(0\) to \(T\) as \(x\) ranges over \([0,1]\). Now, let 
$$u_k(t)=(a^{1/2}\psi_k)(x(t)),\qquad t\in[0,T].$$ As shown in~\cite{everitt2005sturm}, we have that \(u_k\) satisfies
\begin{EQA}\label{uk-eq}
\begin{cases}-u_k''+ Qu_k = \lambda_k^{-1}u_k, \quad t\in(0,T),\\
u_k(0)=0,\quad c_1u_k'(T)+c_2u_k(T)=0,\end{cases} &&
\end{EQA} where \(Q(t)=(q-p^{1/4}(p(p^{-1/4})')')(x(t))\) and \(c_1,c_2\) depend on \(a,b\). Recalling that \(q=b^2-(ab)'\) and \(p=a^2\), we can simplify \(Q\) as
\begin{EQA}
&&Q(t) = (b^2-(ab)' + \tfrac14(a')^2+\tfrac12aa'')(x(t)).
\end{EQA}  The values of \(c_1,c_2\) influence the values \(\lambda_k\) for which a solution to~\eqref{uk-eq} exists. However, we will not need the second boundary condition \(c_1u_k'(T)+c_2u_k(T)=0\) in the subsequent analysis of the eigenfunctions. 
\begin{lemma}\label{lma:psi-v}We have
\begin{EQA}\label{psik0bd}
\|\psi_k\|_\infty &\leq& \frac{\|u_k\|_\infty}{a_0^{1/2}},\\
\|\psi_k'\|_\infty &\leq& \frac{\|a'\|_\infty\|u_k\|_\infty + \|u_k'\|_\infty}{a_0^{3/2}}.\label{psik1bd}
\end{EQA} Here, \(\|\cdot\|_\infty\) are norms over appropriate domains, i.e. \([0,1]\) for \(\psi_k\) and \([0,T]\) for \(u_k\). 
\end{lemma}
This follows by simple calculations using that \(\psi_k(x) = u_k(t(x))a(x)^{-1/2}\), and that \(t'=a^{-1}\). 
\begin{lemma}\label{lma:vkbds}
We have \(u_k'(0)\neq0\). Let \(v_k=u_k/u_k'(0)\). Then for all \(k\) large enough, it holds
\begin{EQA}\label{eq:vkbds}
&&\|v_k\|_\infty\leq 2\sqrt{\lambda_k},\quad \|v_k'\|_\infty\leq 2, \quad \|v_k\|_{L^2}\geq c\sqrt{\lambda_k},
\end{EQA} where \(c\) is a constant depending only on the functions \(a\) and \(b\).
\end{lemma} The proof of Lemma~\ref{lma:vkbds} (below the proof of Theorem~\ref{thm:psik}) is an adaptation from Chapter 2, Section 3.26 of~\cite{yoshida1991lectures}.
\begin{proof}[Proof of Theorem~\ref{thm:psik}]
We have already shown that \(\lambda_k\) have the specified decay rate, in~\eqref{lambdak-asymp}. To bound \(\|\psi_k\|_\infty,\|\psi_k'\|_\infty\), we will write the bounds from Lemma~\ref{lma:psi-v} in terms of \(v_k\), and then apply~\eqref{eq:vkbds}. Note that \(u_k\) has unit \(L^2\) norm. Indeed. 
\begin{EQA}
&&\int_0^Tu_k(t)^2 = \int_0^Ta(x(t))\psi_k(x(t))^2dt = \int_0^1\psi_k(x)^2dx =1,
\end{EQA} using that \(dt=t'(x)dx = a(x)^{-1}dx\). That \(\|u_k\|_{L^2}=1\) implies \(u_k\) can be expressed in terms of \(v_k\) as \(u_k = v_k/\|v_k\|_{L^2}\). Thus the bounds from Lemma~\ref{lma:psi-v} can be written as
\begin{EQA}\label{psik0bd}
\|\psi_k\|_\infty &\leq& \frac{\|v_k\|_\infty}{a_0^{1/2}\|v_k\|_{L^2}},\\
\|\psi_k'\|_\infty &\leq& \frac{\|a'\|_\infty\|v_k\|_\infty + \|v_k'\|_\infty}{a_0^{3/2}\|v_k\|_{L^2}}.\label{psik1bd}
\end{EQA} 
Thus~\eqref{eq:vkbds} now gives \(\|\psi_k\|_\infty \leq C\) and \(\|\psi_k'\|_\infty \leq C\lambda_k^{-1/2}\leq Ck\), using~\eqref{lambdak-asymp}. This holds for all \(k\) large enough. But since \(\psi_k\in C^2[0,1]\), we can extend these two bounds to hold for all \(k\), by increasing the constant \(C\).
\end{proof}

\begin{proof}[Proof of Lemma~\ref{lma:vkbds}]
By standard theory of uniqueness for linear ODEs, the only function \(u\) satisfying \(-u_k''+ Qu_k = \lambda_k^{-1}u_k\) and \(u_k(0)=u_k'(0)=0\) is \(u_k\equiv0\). But since we know \(\psi_k\) is not identically zero, \(u_k\) is also not. Thus we can't have \(u_k'(0)=0\). 

Define \(\rho_k = \lambda_k^{-1/2}\), so that \(v_k''+\rho_k^2v_k = Qv_k\) (since \(u_k\) satisfies this same equation). Now, define the function 
\begin{EQA}
&&\tilde v(t) =  \rho_k^{-1}\int_0^t\sin(\rho_kt-\rho_ks)Q(s)v_k(s)ds.
\end{EQA} It is straightforward to show that \(\tilde v\) satisfies \(\tilde v'' +\rho_k^2\tilde v = Qv_k\). Thus \((v_k-\tilde v)'' =-\rho_k^2(v_k-\tilde v)\), and therefore \((v_k-\tilde v)(t)=A\cos(\rho_kt)+B\sin(\rho_kt)\) for some \(A\) and \(B\). We have therefore shown
\begin{EQA}
&&v_k(t) = A\cos(\rho_kt) + B\sin(\rho_kt) + \rho_k^{-1}\int_0^t\sin(\rho_kt-\rho_ks)Q(s)v_k(s)ds.
\end{EQA} Furthermore, the initial conditions \(v_k(0)=0,v_k'(0)=1\) imply \(A=0\) and \(B=1/\rho_k\). Thus \(v_k\) satisfies
\begin{EQA}\label{vk-fixedpt}
&&v_k(t)=\rho_k^{-1}\sin(\rho_kt) +  \rho_k^{-1}\int_0^t\sin(\rho_kt-\rho_ks)Q(s)v_k(s)ds,
\end{EQA}
which implies
\begin{EQA}\label{vkprime-fixedpt}
&&v_k'(t)=\cos(\rho_kt) + \int_0^t\cos(\rho_kt-\rho_ks)Q(s)v_k(s)ds.
\end{EQA} We now use~\eqref{vk-fixedpt} to bound \(\|v_k\|_\infty\). Specifically, we have
\begin{EQA}
&&\|v_k\|_\infty=\sup_{t\in[0,T]}|v_k(t)| \leq \rho_k^{-1} + \rho_k^{-1}T\|Q\|_\infty\|v_k\|_\infty.
\end{EQA} Using that \(\rho_k=\lambda_k^{-1/2}>2T\|Q\|_\infty\) for all \(k\) large enough, we conclude that
\begin{EQA}\label{vkinfty-bd}
&&\|v_k\|_\infty\leq \frac{\rho_k^{-1}}{1-T\|Q\|_\infty/\rho_k} \leq 2\rho_k^{-1} = 2\sqrt{\lambda_k}
\end{EQA} We now use~\eqref{vkinfty-bd} and~\eqref{vkprime-fixedpt} to obtain a bound on \(\|v_k'\|_\infty\). We have
\begin{EQA}
&&\|v_k'\|_\infty=\sup_{t\in[0,T]}|v_k'(t)| \leq 1 + T\|Q\|_\infty\|v_k\|_\infty \leq 1 + 2T\|Q\|_\infty/\rho_k \leq 2.
\end{EQA} Finally, we obtain a lower bound on \(\|v_k\|_{L^2}\). To do so, consider the ansatz \(v_k(t) = \rho_k^{-1}\sin(\rho_kt) + \rho_k^{-2}w_k(t)\). Plugging this into~\eqref{vk-fixedpt} gives
\begin{EQA}
&&w_k(t) = \int_0^t\sin(\rho_kt-\rho_ks)Q(s)\sin(\rho_ks)ds + \rho_k^{-1}\int_0^t\sin(\rho_kt-\rho_ks)Q(s)w_k(s)ds.
\end{EQA} Thus \(\|w_k\|_\infty \leq T\|Q\|_\infty + \rho_k^{-1}T\|Q\|_\infty\|w_k\|_\infty\), from which we deduce
\begin{EQA}
&&\|w_k\|_\infty \leq\frac{T\|Q\|_\infty}{1-T\|Q\|_\infty/\rho_k} \leq 2T\|Q\|_\infty.
\end{EQA} Now,
\begin{EQA}
\|v_k\|_{L^2} &=& \rho_k^{-1}\left\|\sin(\rho_k\cdot) + \rho_k^{-1}w_k\right\|_{L^2}\geq \rho_k^{-1}\|\sin(\rho_k\cdot)\|_{L^2} - \rho_k^{-2}\|w_k\|_{L^2}
\end{EQA} Next, we have \(\|w_k\|_{L^2}\leq T\|w_k\|_{\infty} \leq 2T^2\|Q\|_\infty\) and a short calculation gives \(\|\sin(\rho_k\cdot)\|_{L^2}^2 =T/2 - \sin(2T\rho_k)/(4\rho_k) \geq T/2-1/(4\rho_k)\). Combining these estimates we obtain
\begin{EQA}
&&\|v_k\|_{L^2} \geq \frac{1}{\rho_k}\sqrt{\frac{T}{2} - \frac{1}{4\rho_k}} - \frac{2T^2\|Q\|_\infty}{\rho_k^2},
\end{EQA} and this quantity can be further bounded below by \(c/\rho_k\) for some \(c>0\) depending only on \(T\) and \(\|Q\|_\infty\) (and hence only on the functions \(a\) and \(b\)).
\end{proof}

\begin{proof}[Proof of Lemma~\ref{lma:RRT}]
Note that for a smooth function \(g\in C^\infty[0,1]\) with \(g(0)=0\), we have \(\calR(ag'+bg)=g\). This immediately follows from the characterization~\eqref{calR-def}, with \(f=ag'+bg\). We claim that \(\calR^\T h\) satisfies \(-(a\calR^\T h)' + b\calR^\T h = h\). To see this, it suffices to show \(\int\left[-(a\calR^\T h)' + b\calR^\T h\right]g = \int hg\) for all \(g\in C_c^\infty[0,1]\). This is straightforward to show using integration by parts, the definition of \(\calR^\T\), and the observation that \(\calR(ag'+bg)=g\). 
Note also that \(\calR^\T h(1)=0\). To see this, take any smooth function \(g\in C^\infty[0,1]\) such that \(g(0)=0,g(1)=1\). Using analogous calculations to those just discussed, we have
\begin{EQA}
\int hg &=& \int\left[-(a\calR^\T h)' + b\calR^\T h\right]g\\
& =& (a\calR^\T h)(0)g(0)-(a\calR^\T h)(1)g(1) + \int h\calR(ag'+bg) \\
& =& -a(1)(\calR^\T h)(1) +\int hg.
\end{EQA} Since \(a(1)\neq0\), we conclude that \((\calR^\T h)(1)=0\). To summarize, 
\begin{EQA}\label{calRT-def}
&&g=\calR^\T h \quad \Leftrightarrow\quad -(ag)' + bg = h,\quad g(1)=0,
\end{EQA}
Finally, we study the operator \(\calR\calR^\T\). Fix \(f\in L^2[0,1]\) and let \(g=\calR^\T f\)  and \(h=\calR g=\calR\calR^\T f\). We first determine the ODE satisfied by \(h\). Using~\eqref{calR-def} and~\eqref{calRT-def}, we have
\begin{EQA}
&&ah'+bh=g,\quad -(ag)'+bg = f,\\
&&\implies -(a[ah'+bh])' + b[ah'+bh]=f.
\end{EQA}
Expanding and rearranging terms, we get \(-(a^2h')' + (b^2-(ab)')h = f\). Next, we determine the boundary conditions \(h\) satisfies. Since \(h=\calR g\), we immediately conclude \(h(0)=0\), by~\eqref{calR-def}. Furthermore, again by~\eqref{calR-def}, we have \(ah'+bh=g\). But since \(g=\calR^\T f\), we have \(g(1)=0\) by~\eqref{calRT-def}. Thus \(a(1)h'(1)+b(1)h(1)=g(1)=0\). Combining all of the above, we conclude that \(\calR\calR^\T\) is as in~\eqref{calRRT-def}.
\end{proof}

\section{Auxiliary results}\label{app:aux}


Let \(a_k,b_k,c_k>0\), \(k=1,\dots,\dimp\) with \(b_k\geq c_k\) for \(k=1,\dots,m\) and \(b_k\leq c_k\) for \(k=m+1,\dots,\dimp\). Then
\begin{EQA}\label{splitsum}
\frac12\sum_{k=1}^{m}\frac{a_k}{b_k}+\frac12\sum_{k=m+1}^{\dimp}\frac{a_k}{c_k}&\leq&\sum_{k=1}^{\dimp}\frac{a_k}{b_k+c_k}\leq \sum_{k=1}^{m}\frac{a_k}{b_k}+\sum_{k=m+1}^{\dimp}\frac{a_k}{c_k}.
\end{EQA}
Let\(b> a\geq1\) be whole numbers. If \(\alpha>-1\) then 
\begin{EQA}\label{asympsum1}
\frac{b^{\alpha+1}-a^{\alpha+1}}{\alpha+1}&\leq&\sum_{k=a}^bk^\alpha \leq \frac{(b+1)^{\alpha+1}}{\alpha+1}.
\end{EQA}
If \(\alpha<-1\) then
\begin{EQA}\label{ab-alpha-1}
1\leq \sum_{k=1}^bk^{\alpha}&\leq& 1 + \frac{1}{-1-\alpha}.
\end{EQA}
If \(\alpha<-1\) and \(a\geq1\), then
\begin{EQA}\label{tail-sum}
\frac{(a+1)^{\alpha+1}-(b+1)^{\alpha+1}}{-1-\alpha} &\leq &\sum_{k=a+1}^bk^\alpha \leq \frac{a^{\alpha+1}}{-1-\alpha}
\end{EQA} 
\begin{lemma}\label{lma:AB}Suppose \(0\preceq A_1\preceq A \preceq A_2\) and  \(0\preceq B_1\preceq B \preceq B_2\). Then
\begin{EQA}
\|A_2^{-1/2}B_1A_2^{-1/2}\|&\leq&\|A^{-1/2}BA^{-1/2}\|\leq \|A_1^{-1/2}B_2A_1^{-1/2}\|,\\
\Tr(A_2^{-1/2}B_1A_2^{-1/2})&\leq &\Tr(A^{-1/2}BA^{-1/2})\leq\Tr(A_1^{-1/2}B_2A_1^{-1/2}).
\end{EQA}
\end{lemma}
\begin{proof}
We have
\begin{EQA}
\|A^{-1/2}BA^{-1/2}\|&=&\sup_{\xv\neq0}\frac{\xv^\T A^{-1/2}BA^{-1/2}\xv}{\xv^\T\xv}=\sup_{\yv\neq0}\frac{\yv^\T B\yv}{\yv^\T A\yv}
\end{EQA} Using the matrix inequalities now gives
\begin{EQA}
\frac{\yv^\T B_1\yv}{\yv^\T A_2\yv}&\leq&\frac{\yv^\T B\yv}{\yv^\T A\yv}\leq\frac{\yv^\T B_2\yv}{\yv^\T A_1\yv}
\end{EQA} Taking the supremum over all \(\yv\neq0\) proves the operator norm inequalities. For the trace inequalities, first note that \(A^{-1/2}B_1A^{-1/2}\preceq A^{-1/2}BA^{-1/2}\preceq A^{-1/2}B_2A^{-1/2}\). Therefore, 
\begin{EQA}
\Tr(A^{-1/2}B_1A^{-1/2})&\leq &\Tr(A^{-1/2}BA^{-1/2})\leq\Tr(A^{-1/2}B_2A^{-1/2}).
\end{EQA} Next, we have 
\begin{EQA}
\Tr(A^{-1/2}B_2A^{-1/2})&=&\Tr(B_2^{1/2}A^{-1}B_2^{1/2})\leq\Tr(B_2^{1/2}A_1^{-1}B_2^{1/2})=\Tr(A_1^{-1/2}B_2A_1^{-1/2}).
\end{EQA}
To get the inequality, we used that \(B_2^{1/2}A^{-1}B_2^{1/2}\preceq B_2^{1/2}A_1^{-1}B_2^{1/2}\). This proves the desired righthand trace inequality. The lefthand trace inequality is proved analogously.
\end{proof}

The following result is a concentration inequality for a Gaussian vector, from~\cite{laurent2000adaptive}. 
\begin{theorem}\label{TexpbLGA}
Let \( \gaussv \sim \ND(0,\Id_{\dimp}) \) and let \(\BBH\succ0\) satisfy \(\|\BBH\|=1\), \( \dimH = \tr(\BBH) \), and \( \vH^{2} = \tr(\BBH^{2}) \). Then for each \( \xx \geq 0 \), it holds
\begin{EQA}
\label{Pxiv2dimAvp12}
	\P\Bigl( \gaussv^{\T}\BBH\gaussv - \dimH > 2 \vH \, \sqrt{\xx} + 2 \xx \Bigr)
	& \leq &
	\ex^{-\xx} .
\label{Pxiv2dimAvp12m}
\end{EQA}
\end{theorem}

Since \(\|\BBH\|=1\), it follows that \(\vH^{2}=\Tr(\BBH^2) \leq\Tr(\BBH)=\dimH\). Using this inequality in Theorem~\ref{TexpbLGA} yields the following corollary.
\begin{corollary}\label{corr:gaussconc}
Let \(\BBH\succ0\) and \( \|\BBH\|=1\). Define \(\dimH = \Tr(\BBH)\). Then
\begin{EQA}
\label{gaussconc}
	\P\Bigl( \gaussv^{\T}\BBH\gaussv  > (\sqrt{\dimH}+t)^2\Bigr)
	& \leq &
	\ex^{-t^2/2} .
\end{EQA}
\end{corollary}
\begin{proof}
Let \(t=\sqrt{2\xx}\) and recall \( \vH^{2} = \tr(\BBH^{2}) \). Using that \(\vH^2\leq\dimH\),  we have
\begin{EQA}
\P\Bigl( \gaussv^{\T}\BBH\gaussv  > (\sqrt{\dimH}+t)^2\Bigr)&=& \P\Bigl( \gaussv^{\T}\BBH\gaussv  -\dimH > 2\sqrt{2\xx}\sqrt{\dimH} +2\xx\Bigr)\\
&\leq &\P\Bigl( \gaussv^{\T}\BBH\gaussv  -\dimH>2\vH\sqrt{\xx}\ + 2\xx\Bigr) \leq \ex^{-\xx}=\ex^{-t^2/2}.
\end{EQA} 
\end{proof}

\subsection{Taylor remainder inequalities}\label{app:subsec:taylor}
Recall the following three definitions:
\begin{EQA}
\dltwb(\UV,\DFL) &=&\sup_{\thv\in\UV}\frac{|\lgd(\thv)-\lgd(\map)-\|\DFLG(\thv-\map)\|^2/2|}{\|\DFL(\thv-\map)\|^2/2},\label{def:omega:app}\\
\dltwb_3(\UV,\DFL) &=&  \sup_{\thv\in\UV}\bigl\|\DFL^{-1}\bigl(\nabla^2\lgd(\thv)-\nabla^2\lgd(\map)\bigl)\DFL^{-1}\bigr\|,\label{def:omega3:app}\\
\dltwu_3(\UV,\DFL) &=& \sup_{\thv\in\UV}\bigl\|\DFL^{-1}\bigl(\nabla^2\lgd(\thv)-\nabla^2\lgd(\map)\bigl)\DFL^{-1}\bigr\|/\|\DFL(\thv-\map)\|.\label{def:tau3:app}
\end{EQA}
We claim that 
\begin{EQA}
\dltwb(\UV,\DFL)& \leq& \frac13\sup_{\thv\in\UV}\|\DFL(\thv-\map)\|\dltwu_3(\UV,\DFL),\label{del31}\\
\dltwu_3(\UV,\DFL)&\leq&\sup_{\thv\in\UV}\|\nabla^3\lgd(\thv)\|_{\DFL},\label{del32}\\
\dltwb_3(\UV,\DFL)&\leq &\sup_{\thv\in\UV}\|\DFL(\thv-\map)\|\dltwu_3(\UV,\DFL).\label{del33}
\end{EQA}
~\eqref{del31} and~\eqref{del32} prove~\eqref{del3bds1}, since \(\sup_{\thv\in\UV(\DFc,\rrdc)}\|\DFc(\thv-\map)\|=\rrdc\). Also,~\eqref{del32} and~\eqref{del33} prove~\eqref{del3bds2}. To prove~\eqref{del31}, fix \(\thv\in\UV\), let \(\uv=\thv-\map\), and define \(h(t) = \lgd(\map+t\uv)\). Then using the integral form of the Taylor remainder theorem, we have
\begin{EQA}
|\lgd(\thv)-\lgd(\map)&-&\|\DFLG(\thv-\map)\|^2/2 = (h(1)-h(0)-h'(0))-\frac12h''(0) \\
&=&\int_0^1(1-t)h''(t)dt - \frac12h''(0) = \int_0^1(1-t)(h''(t)-h''(0))dt\\
&=&\int_0^1(1-t)\bigl\langle\nabla^2\lgd(\map+t\uv)-\nabla^2\lgd(\map), \uv^{\otimes2}\bigr\rangle dt.
\end{EQA}
Thus
\begin{EQA}
|\lgd(\thv)-\lgd(\map)&-&\|\DFLG(\thv-\map)\|^2/2  \\
&\leq &\int_0^1(1-t)\|\DFL^{-1}(\nabla^2\lgd(\map+t\uv)-\nabla^2\lgd(\map))\DFL^{-1}\|\|\DFL\uv\|^2dt\\
&\leq &\dltwu_3(\UV, \DFL)\|\DFL\uv\|^3\int_0^1(1-t)tdt = \dltwu_3(\UV,\DFL)\|\DFL\uv\|^3/6
\end{EQA} for all \(\uv=\thv-\map\) where \(\thv\in\UV\). We conclude~\eqref{del31} using the definition of \(\dltwb(\UV,\DFL)\). To prove~\eqref{del32}, note that
\begin{EQA}
&&\left\|\DFL^{-1}\left(\nabla^2\lgd(\thv)-\nabla^2\lgd(\map)\right)\DFL^{-1}\right\|=\sup_{\|\DFL\vv\|=\|\DFL\wv\|=1}\langle \nabla^2\lgd(\thv)-\nabla^2\lgd(\map),\;\vv\otimes\wv\rangle.
\end{EQA} Thus
\begin{EQA}\label{dltwsupuv}
&&\dltwu_3(\UV,\DFL) = \sup_{\thv\in\UV}\sup_{\|\DFL\vv\|=\|\DFL\wv\|=1}\frac{\langle  \nabla^2\lgd(\thv)-\nabla^2\lgd(\map),\;\vv\otimes\wv\rangle}{\|\DFL(\thv-\map)\|}
\end{EQA} Now, let \(h(t)=\langle\nabla^2\lgd(\map+t\uv), \vv\otimes\wv\rangle\), where \(\uv=\thv-\map\). Then 
\begin{EQA}\label{biglangleA}
\bigl|\langle  \nabla^2\lgd(\map+\uv)&-&\nabla^2\lgd(\map),\;\vv\otimes\wv\rangle\bigr| = |h(1)-h(0)|\leq\sup_{t\in[0,1]}|h'(t)|\\
& =& \sup_{t\in [0,1]}\bigl|\langle \nabla^3\lgd(\map+t\uv),\uv\otimes\vv\otimes\wv\rangle\bigr|\leq\|\nabla^3\lgd(\map+t\uv)\|_{\DFL}\|\DFL\uv\|,
\end{EQA} using the definition~\eqref{TDdef} of \(\|T\|_{\DFL}\) for a tensor \(T\), and the fact that \(\|\DFL\vv\|=\|\DFL\wv\|=1\). Substituting~\eqref{biglangleA} into~\eqref{dltwsupuv} concludes the proof of~\eqref{del32}. Finally,~\eqref{del33} is immediate from the definition of \(\dltwu_3\).

\bibliography{bibliogr_Lap}  

\end{document}